\documentclass{birkjour}

\usepackage[noadjust]{cite}
\usepackage{xcolor}
\RequirePackage[all]{xy}

\newcommand{\BibTeX}{B\kern-0.1emi\kern-0.017emb\kern-0.15em\TeX}
\newcommand{\XYpic}{$\mathrm{X\kern-0.3em\raisebox{-0.18em}{Y}}$-$\mathrm{pic}\,$}

\usepackage{amsbsy,amsfonts,amsmath,amssymb,enumerate,epsfig,graphicx,rotating,subfigure}

\usepackage{amsmath}
\usepackage{amsfonts}
\usepackage{amssymb}
\usepackage{amsthm}

\numberwithin{equation}{section}

\newtheorem{theorem}{Theorem}[section]

\newtheorem{corollary}[theorem]{Corollary}
\newtheorem{proposition}[theorem]{Proposition}

\newtheorem{remark}{Remark}

\begin{document}

\title[On a class of Rainville type generating functions for classical OP]{On a class of Rainville type generating functions for classical orthogonal polynomials}
\author[M. B. Zahaf ]{Mohammed Brahim Zahaf    }
\address{Department of mathematics\\ Abou Bekr Belkaid University \\ Tlemcen 13000\\ Algeria\\
and\\ Laboratory of nonlinear analysis   and applied  mathematics\\ Abou Bekr Belkaid University\\
	Tlemcen 13000 \\ Algeria}
\email{m\_b\_zahaf@yahoo.fr}

	 \author[M. Mesk]{Mohammed Mesk}
\address{Department of ecology and environment\\
	Abou Bekr Belkaid University\\  Tlemcen 13000\\ Algeria\\
and\\ Laboratory of nonlinear analysis   and applied  mathematics\\ Abou Bekr Belkaid University\\
	Tlemcen 13000 \\ Algeria}
\email{m\_mesk@yahoo.fr }



\date{\today}
 
\begin{abstract}
In this paper, we used algebraic calculus to characterize all generating functions of the form $A(t)F(xtA(t)-R(t))$ for the classical orthogonal polynomials. We also applied the derivative operator to get generating functions of the form $A(t)^kG(xtA(t)-R(t))$ for these polynomials. Particularly, we obtained bilateral generating functions between Hermite and associate Hermite polynomials and between ultraspherical and Chebyshev polynomials of the second kind.  
\end{abstract}

\maketitle


\section{Introduction}
For a natural number $n$ and a real $x$, we say that $\{P_n(x)\}$ is a monic polynomial set (MPS) if $P_n(x)=x^n+\text{lower terms}$. A generating function (GF) of Rainville type for a MPS $\{P_n(x)\}$ has the form \cite{rainville1971special}
\begin{equation}\label{gf1}
B(t)F\left(xtA(t)-R(t)\right)=\sum_{n\geq 0}\alpha_nP_n(x)t^n,\,\,(\alpha_n\neq 0 \text{ for } n\geq 0),
\end{equation}
where $F(t)=\sum_{n\geq 0}\alpha_nt^n$, $A(t)$, $B(t)$ and $R(t)$ are formal power series with  $(F(0),A(0),B(0),R(0))=(1,1,1,0)$. The GF \eqref{gf1}  is  a generalization of the Boas-Buck type GF $B(t)F(xtA(t))$ \cite{Boas_Buck_1}. These two forms of GFs generate various interesting polynomial sets, among them, MPSs satisfying the $2$-order recursion
\begin{eqnarray}\label{gf6}
	\left\{
	\begin{array}{l}
		xP_n(x)=P_{n+1}(x)+\beta_n P_n(x)+\omega_n P_{n-1}(x),\quad n\geq
		0,\\
		P_{-1}(x)=0,\;\;P_0(x)=1
	\end{array} 
	\right.
\end{eqnarray} 
where $\{\beta_n\}$ and $\{\omega_n\}$ are complex sequences. Note that with some additional conditions the recursion \eqref{gf6} characterizes orthogonal polynomials, see  \cite[p.21,Theorem 4.4]{Chihara}.  A popular subclass is the classical orthogonal polynomials  which have the property that their derivatives are   also orthogonal \cite{sonine1887uber}, namely, Jacobi,  Gegenbauer, Chebyshev of the first and second kinds, Hermite, Laguerre and Bessel polynomials \cite{castillo2024first,ismail2005classical,srivastava2023introductory}. Orthogonal polynomials can be found in a wide range of fields, including applied mathematics, mathematical physics, and approximation theory \cite{nikiforov1988special, raslan2022operational, xu2021approximation}.

The problem of characterizing all orthogonal polynomials with a Boas-Buck or Rainville type is not yet completely solved. In the meanwhile, some specific situations have been investigated by many authors, see for instance \cite{alsalam,ChiharaBrenke68,ismail1974obtaining,Kubo13generalorder,meskzahaf4,meskzahaf,Meixner1934,Sheffer37,Bachhaus}. We also quote further extensions to higher order recursions, see \cite{Bencheikh,BencheikhGam16,CHAGGARA2024128123,MeskZahaf_d_Symmetric,meskzahaf_1} and the references therein. In \cite{ismail1974obtaining} Ismail gave a method for obtaining GFs of Boas Buck type for orthogonal polynomials and applied it to the classical polynomials. In \cite{meskzahaf4} we gave a method for obtaining MPSs satisfying \eqref{gf6} and generated by the Rainville type GF \eqref{gf1} with $B(t)=A(t)$,i.e, 
\begin{equation}\label{gf0}
	A(t)F\left(xtA(t)-R(t)\right)=\sum_{n\geq 0}\alpha_nP_n(x)t^n.
\end{equation}
Motivated by Ismail's work \cite{ismail1974obtaining}, we use the method stated in \cite{meskzahaf4} to classify all generating functions \eqref{gf0} of the classical orthogonal polynomials. We also exploit the fact that the derivative of a classical orthogonal polynomial is also classical to obtain some new and known Rainville type GFs of the form \eqref{gf1}.

In section \ref{sec2}, we give some  preliminary results where some of them were  established in \cite{meskzahaf4} as the relations between $\{\beta_n\}$, $\{\omega_{n}\}$ and the coefficients of the formal series of powers $A(t)$, $F(t)$ and $R(t)$. In section \ref{sec3} we apply the results of section \ref{sec2} to classify the GFs \eqref{gf0} of the classical orthogonal polynomials. In section \ref{sec4} we give some GFs of Rainville type \eqref{gf1}  by using the derivative operator.

\section{ Preliminary results}\label{sec2}
\begin{proposition}\cite{meskzahaf4}\label{prop1}
	Let $\{P_n\}_{n\geq 0}$ be a MPS generated by \eqref{gf0},  
	$\frac{A'(t)}{A(t)}=\sum_{n\geq 0}\ S_{n}t^{n}$ and $\frac{R'(t)}{A(t)}=\sum_{n\geq 0}\ R_{n}t^{n}$. Then we have 
	\begin{eqnarray}\label{gf7}
	\alpha_nxP'_n(x)-n\alpha_{n}P_{n}(x)&=&-\sum_{k=0}^{n-1}S_{k}\alpha_{n-k-1}\left[xP'_{n-k-1}(x)+P_{n-k-1}(x)\right]\nonumber \\
	&&+\sum_{k=0}^{n-1}R_{k}\alpha_{n-k}P'_{n-k}(x),\;\;
	n\geq 1.
	\end{eqnarray}
\end{proposition}

\begin{proposition}\cite{meskzahaf4}\label{prop2}
	For the MPS generated by \eqref{gf0} and satisfying \eqref{gf6} we have:
	\begin{equation}\label{gf711}
	xP'_n(x)-nP_{n}(x)=\sum_{k=0}^{n-1}A^k_nP'_{n-k}(x),\;\;
	n\geq 1,
	\end{equation}
	where
	\begin{eqnarray}
	A^0_n&=&-S_0\frac{\alpha_{n-1}}{\alpha_n}+R_0, \;n\geq1,\\
	A^1_n&=&-S_1\frac{\alpha_{n-2}}{\alpha_n}-S_0\frac{\alpha_{n-1}}{\alpha_n}\beta_{n-1}+R_{1}\frac{\alpha_{n-1}}{\alpha_n},\; n\geq 2,\label{A1n11}\\ 
A^k_n&=&-S_{k}\frac{\alpha_{n-k-1}}{\alpha_n}-S_{k-1}\frac{\alpha_{n-k}}{\alpha_n}\beta_{n-k}
-S_{k-2}\frac{\alpha_{n-k+1}}{\alpha_n}\omega_{n-k+1}\nonumber\\
&&+R_{k}\frac{\alpha_{n-k}}{\alpha_n},\;\;  k\geq 2, \,n\geq k+1.\label{Aknnn}
	\end{eqnarray}

\end{proposition}

\begin{proposition}\cite{meskzahaf4}\label{prop3}
	For a MPS $\{P_n\}$ satisfying \eqref{gf6} and \eqref{gf0} with $(n-k)A^{k}_{n}=0$ and $n=k\geq 0$, we have 
	
	i.)
	\begin{eqnarray}\label{eqi}
	\beta_n &=&(n+1)A^{0}_{n+1}-nA^{0}_{n},\;n\geq 0.
	\end{eqnarray}
	
	ii.)
	\begin{eqnarray}\label{eqii}
	\,\omega_n &=& \frac{n}{2}A^1_{n+1}-\frac{n-1}{2}A^1_n-\frac{n}{2(n+1)}\left(\beta_n-A^0_n\right)^2,\;n\geq 1.
	\end{eqnarray}
	
	iii.)
	\begin{eqnarray}\label{eqiii}
	&&A^2_{n+1}-\frac{n-2}{n-1}A^2_n=\omega_n\left(\frac{1}{n+1}\beta_n+\frac{2}{n}\beta_{n-1}-\frac{n+2}{n}A^0_{n-1}+\frac{n}{n+1}A^0_n\right)\nonumber
	\\
	&&+A^1_n\left(-\frac{n+2}{n+1}\beta_n+\frac{n-1}{n}\beta_{n-1}+\frac{1}{n+1}A^0_n+\frac{1}{n}A^0_{n-1}\right),\;n\geq 2. 
	\end{eqnarray}
	
	iv.)
	\begin{eqnarray}\label{eqiv}
	&&\frac{2}{n}\omega_n\omega_{n-1} = A^3_{n+1}-\frac{n-3}{n-2}A^3_n+
	\frac{n+2}{n}\omega_n A^1_{n-1}-\frac{n-1}{n}\omega_{n-1}A^1_n-\frac{1}{n}A^1_nA^1_{n-1}\nonumber\\
	&&-\left(-\frac{n+2}{n+1}\beta_n +\frac{n-2}{n-1}\beta_{n-2}+\frac{1}{n+1}A^0_n+\frac{1}{n-1}A^0_{n-2}\right)A^2_n,\;n\geq 3.
	\end{eqnarray}
	
	v.)
	\begin{eqnarray}\label{eqv}
	&&A^k_{n+1}-\frac{n-k}{n-k+1}A^k_n+\left(\frac{n+2}{n+1}\beta_n-\frac{n-k+1}{n-k+2}\beta_{n-k+1}\right)A^{k-1}_n\nonumber\\
	&&\quad\qquad\qquad+\frac{n+2}{n}\,\omega_n\, A^{k-2}_{n-1}-\frac{n-k+2}{n-k+3}\,\omega_{n-k+2}\,A^{k-2}_n\nonumber\\
	&&\qquad\qquad\qquad=\sum_{l=0}^{k-1}\frac{A^{k-l-1}_n\,A^l_{n-k+l+1}}{n-k+l+2},\;\; n\geq k\geq 4.
	\end{eqnarray}

\end{proposition}

\begin{corollary}\label{corl1}
If $A_n^2=0$ for $n\geq 3$ and $A_n^3=0$ for $n\geq 4$ then $A_n^{k}=0$ for $k\geq 4 $, $n\geq k+1$.
\end{corollary}
 \begin{proof}
 We will use \eqref{eqv}  and proceed by induction on $k$ to show that $A_n^{k}=0$ for $k\geq4 $, $n\geq k+1$. Indeed $k=4$ in \eqref{eqv} leads to 
\begin{equation}\label{AAS4}
A^4_{n+1}-\frac{n-4}{n-3}A^4_n=0, \;n\geq4.
\end{equation}
By summing \eqref{AAS4} after multiplying by $n-3$ we obtain $A_n^4=0$ for $n\geq 5$. Suppose that $A_n^4=A_n^5=...=A_n^{k}=0$ and  by using the shifting $k\to k+1$ in \eqref{eqv} we obtain 
\begin{equation}\label{AAkS}
A^{k+1}_{n+1}-\frac{n-k-1}{n-k}A^{k+1}_n=0, \;n\geq k+1.
\end{equation}
By summing \eqref{AAkS} after multiplying by $n-k$ we obtain $A_n^{k+1}=0$ for $n\geq k+2$.

 \end{proof}
\begin{corollary}\label{cor1}
If  $\{P_n\}$ is a symmetric MPS generated by \eqref{gf0} and satisfying  \eqref{gf6} then we have

i.)
	\begin{equation}\label{eqis}
	A^{2k}_{n}=0,\;k\geq 0,\,n\geq 2k+1.
	\end{equation}
	
	ii.)
	\begin{equation}\label{eqiis}
	\,\omega_n = \frac{n}{2}A^1_{n+1}-\frac{n-1}{2}A^1_n,\;n\geq 1.
	\end{equation}
	
	iii.)
	\begin{equation}\label{eqivs}
	\frac{2}{n}\omega_n\omega_{n-1} =  A^3_{n+1}-\frac{n-3}{n-2}A^3_n+
	\frac{n+2}{n}\omega_n A^1_{n-1}-\frac{n-1}{n}\omega_{n-1}A^1_n-\frac{1}{n}A^1_nA^1_{n-1},\;n\geq 3.
	\end{equation}
	
	iv.)
	\begin{eqnarray}\label{eqvs}
	&&A^{2k+1}_{n+1}-\frac{n-2k-1}{n-2k}A^{2k+1}_n+\frac{n+2}{n}\,\omega_n\, A^{2k-1}_{n-1}-\frac{n-2k+1}{n-2k+2}\,\omega_{n-2k+1}\,A^{2k-1}_n
	\nonumber\\
	&&=\sum_{l=0}^{k-1}\frac{A^{2k-2l-1}_n\,A^{2l+1}_{n-2k+2l+1}}{n-2k+2l+2},\; n\geq 2k+1,\,k\geq 2.
	\end{eqnarray}

\end{corollary}
\begin{proof}
$\{P_n\}$ is a symmetric PS means that  $P_n(-x)=(-1)^nP_n(x)$ for $n\geq 0$. The substitution $x\rightarrow -x$ in equation \eqref{gf711} minus equation \eqref{gf711} itself left us with
$$(1-(-1)^{k+1})A^k_n=0,\;\text{ for  } 0 \leq k \leq n - 1.$$
So, we have $A^{2l}_{n}=0,\;l\geq 0,\,n\geq 2l+1.$\\
Finally, by using  $A_n^{2k}=0$ in \eqref{eqi}, \eqref{eqii}, \eqref{eqiv} and \eqref{eqv} we obtain \eqref{eqiis}, \eqref{eqivs} and \eqref{eqvs}.
\end{proof}
\begin{remark} If $\{P_n\}$ is as in corollary \eqref{cor1} then we have $S_{2k}=R_{2k}=0$ for  $k\geq 0$.
\end{remark}
\section{Obtaining generating functions for the classical orthogonal polynomials}\label{sec3}
 We begin by defining the generalized hypergeometric function ${}_{p}F_{q}$, where $p,q$ are natural numbers. It has the form:
\begin{eqnarray*}\label{F01}
{}_{p}F_{q}\left(
\begin{array}{llll}
\mu_1,\mu_2,&...,&\mu_p\\
\nu_1,\nu_2,&...,&\nu_q
\end{array};z\right)=\sum_{n\geq 0}\frac{(\mu_1)_n(\mu_2)_n\cdots(\mu_p)_n}{(\nu_1)_n(\nu_2)_n\cdots(\nu_q)_n}\frac{z^n}{n!}\quad\quad
\end{eqnarray*}
where $\mu_1,\mu_{2},...,\mu_p,\nu_1,\nu_{2},...,\nu_q $ are complex parameters.  The symbol $(\mu)_n$ stands for the shifted factorials, i.e., 
$(\mu)_0=1,\;\;(\mu)_n=\mu(\mu+1)\cdots(\mu+n-1),\;\;n\geq 1$.
We also put $a_n=\frac{\alpha_{n}}{\alpha_{n+1}}$ for $n\geq 0$.
\subsection{Ultraspherical polynomials}
The monic Ultraspherical  (or Gegenbauer) polynomials are defined by  (Hypergeometric representation)
\begin{equation*}
C_n^{(\lambda)}(x)=\frac{(2\lambda)_n}{2^n(\lambda)_n}\,{}_2F_1\left(
\begin{array}{cr}
	-n\,\,\,\,n+2\lambda&\\
	&;\frac{1-x}{2} \\
	\lambda+\frac12&
\end{array}\right),\; \lambda\geq -\frac{1}{2}.
\end{equation*}
Notice that the monic Chebyshev polynomials $T_n(x)$ and $U_n(x)$ of the first and second kinds, respectively, are special cases of the monic Ultraspherical  polynomials with $\lambda=0$ and $\lambda=1$ respectively.
\begin{theorem}\label{th1}
	The only generating functions of the form \eqref{gf0} of the monic Gegenbauer polynomials $ C_n^{(\lambda)}(x) $, Chebyshev polynomials of the first and second kinds $T_n(x)$, $U_n(x)$ are
	\begin{eqnarray}	   
	&&\bullet 	\sum_{n=0}^{+\infty}\frac{(\lambda)_n}{n!} C_n^{(\lambda)}(x)t^n=\left(1-tx+\frac{t^2}{4}\right)^{-\lambda},\label{Ultra1}\\
	&& \mbox{with } A(t)=1,\, R(t)=\frac{t^2	}{4} \mbox{ and } F(t)= (1-t)^{-\lambda}.\nonumber \\
	&&\bullet 	\sum_{n=0}^{+\infty}\frac{(\lambda+1)_n}{n!} C_n^{(\lambda)}(x)t^n=\left(1-\frac{t^2}{4}\right)\left(1-tx+\frac{t^2}4\right)^{-\lambda-1},\label{Ultra2}	\\
	&& \mbox{with }A(t)=\left(1-\frac{t^2}{4}\right)^{-\frac{1}{\lambda}},\, R(t)=\left(1+\frac{t^2}{4}\right)A(t)-1  \text{ and }
F( t) =  \left( 1-t \right) 
	^{-\lambda-1}.\nonumber\\
		&&\bullet \sum_{n=1}^{+\infty}\frac{1}{n} T_n(x)t^n=-\ln\left(1-xt+\frac{t^2}{4}\right), \label{T1}\\
	&& \mbox{with } A(t)=1,\, R(t)=\frac{t^2	}{4} \mbox{ and } F(t)=1-\ln(1-t).\nonumber\\
		&&\bullet \sum_{n=0}^{+\infty}U_n(x)t^n=\frac{1}{1-xt+\frac{t^2}{4}},\label{U3}\\
	&&\mbox{with }R(t)= \left( 1+\frac{t^{2}}{4}\right)\,A ( t ) -1, A(t) \mbox{ arbitrary and } F(t)=\frac{1}{1-t}.\nonumber
		\end{eqnarray}
	\begin{eqnarray}
	&& \bullet 1+\alpha_1\sum_{n=1}^{+\infty}2^{n-1}\rho^{n-1}U_n(x)t^n= {\frac { 1+\left(\alpha_1 -2\rho \right) (xt-\frac{\rho}{2}t^2)
		}{1-2\rho xt+{\rho}^{2}{t}^{2}}}\label{U4}\\
	&& \mbox{with } A(t)=1,\, R(t)=\frac{\rho t^2	}{2} \mbox{ and } F(t)=\frac{1+(\alpha_1-2\rho)t}{1-2\rho t}.\nonumber\\
		&&\bullet\sum_{n=0}^{+\infty}\alpha_n U_n(x)t^n={\frac {4-t^{2}+ 2\left( \alpha_{1}-2\,\rho \right)  \left( 2\,xt-\rho\,
				{t}^{2} \right) }{4-8\,\rho\,xt+ \left( 4\,{\rho}^{2}+4\,{x}^{2}
				-2 \right) {t}^{2}-2\,\rho\,x{t}^{3}+\frac{{t}^{4}}{4}}}
		\label{U1}\\
	&&\mbox{with }A(t)=\frac{2}{2-t^2},\, R(t)=\frac{\rho t^2}{2-t^2} \text{ and }	F ( t ) ={\frac { 1+\left(\alpha_1 -2\rho\right) t}{1-2\rho t+2{t}^{2}}}.\nonumber
			\end{eqnarray}
	where $
	\alpha_n=2^n[U_n\left(\rho\right)+(\tfrac{\alpha_1}{2}-\rho)U_{n-1}\left(\rho\right)], \,n\geq 1,\;\alpha_0=1.$\\ 
\end{theorem}
	\begin{remark} From \eqref{U1} we obtain the following GFs 
	\begin{eqnarray}\label{U2}
		&&\sum_{n=0}^{+\infty}2^nU_n\left(\rho\right) U_n(x)t^n={\frac {4-\,{t}^{2} }{4-8\,\rho\,xt+ \left( 4\,{\rho}^{2}+4\,{x}^{2}
				-2 \right) {t}^{2}-2\,\rho\,x{t}^{3}+\tfrac{{t}^{4}}{4}}}\;\;\nonumber\\
			&&	\sum_{n=1}^{+\infty}2^nU_{n-1}\left(\rho\right) U_n(x)t^n={\frac {4\left( 2\,xt-\rho\,
					{t}^{2} \right) }{4-8\,\rho\,xt+ \left( 4\,{\rho}^{2}+4\,{x}^{2}
					-2 \right) {t}^{2}-2\,\rho\,x{t}^{3}+\tfrac{{t}^{4}}{4}}}.\;\;\nonumber
	\end{eqnarray}
\end{remark}

\begin{proof}[Proof of  Theorem~\ref{th1}]
The Jacobi-Szeg\"o parameters of  $C_n^{(\lambda)}(x)$ are given by
  $\beta_n=0$ (symmetric) and 
\begin{equation}\label{omg1}
\omega_n=\frac14\,{\frac {n ( n-1+2\lambda ) }{ ( n+\lambda ) 
 ( n-1+\lambda) }},\;n\geq1.
\end{equation}
Then by  \eqref{eqiis} we obtain
\begin{equation}\label{A1n111111}
A_n^1=\frac{n}{2(n-1+\lambda)}, n\geq 2,
\end{equation}
and equation \eqref{A1n11} can be written as
 \begin{equation}\label{S1eq}
-S_1\frac{\alpha_{n-2}}{\alpha_n}+R_{1}\frac{\alpha_{n-1}}{\alpha_n}=\frac{n}{2(n-1+\lambda)}, n\geq 2.
\end{equation}
By \eqref{A1n111111} and \eqref{eqivs} we get $A_n^3=0$, $n\geq 4$ and according to Corollary \ref{corl1} we have $A_n^{2k+1}=0$ for $k\geq 1$. Then by using \eqref{Aknnn} and multiplying by $\frac{\alpha_n}{\alpha_{n-2k}}$ we find that
\begin{equation*}\label{Skeqk}
-S_{2k+1}\,\frac{\alpha_{n-2k-2}}{\alpha_{n-2k}}-S_{2k-1}\,\omega_{n-2k}\,+R_{2k+1}\,\frac{\alpha_{n-2k-1}}{\alpha_{n-2k}}=0,k\geq 1,\, n\geq 2k+2
\end{equation*}
which by the shifting $n\to n+2k$ becomes
\begin{equation}\label{Skeqk2}
-S_{2k+1}\,\frac{\alpha_{n-2}}{\alpha_{n}}-S_{2k-1}\,\omega_{n}\,+R_{2k+1}\,\frac{\alpha_{n-1}}{\alpha_{n}}=0,k\geq 1,\, n\geq 2.
\end{equation}
Equations \eqref{S1eq} and \eqref{Skeqk2} can be written, respectively, as
\begin{equation}\label{S1eqan}
-S_1\,a_{n-2}\,a_{n-1}+R_{1}\,a_{n-1}=\frac{n}{2(n-1+\lambda)}, n\geq 2,
\end{equation}
and
   \begin{equation}\label{Skeqanshft}
-S_{2k+1}\,a_{n-2}\,a_{n-1}-S_{2k-1}\,\omega_{n}+R_{2k+1}\,a_{n-1}=0,\;k\geq 1,\, n\geq 2.
\end{equation}
Since   $(S_1,R_1)\neq (0,0)$ from \eqref{S1eqan}, we will consider the two cases $S_1\neq 0$ and $S_1=0$.
 
 \textbf{Case 1.  $S_1\neq 0$.} 
 Denote \eqref{Skeqanshft} by $E(n,k)$ and make the operation $a_{1}\,E(n+1,k)-a_n\,E(2,k)$  to obtain
\begin{equation}\label{S2kp1}
-a_{1}a_n(a_{n-1}-a_{0})\,S_{2k+1}+(a_n\omega_2-a_{1}\omega_{n+1})\,S_{2k-1}=0,\;k\geq 1,\, n\geq 2.
\end{equation}
We distinguish the two sub-cases, $\lambda\neq 1$ and $\lambda=1$.  
 
\textbf{ i) } If $\lambda\neq 1$ then, by \eqref{S1eqan}, $a_n$ could not be a constant sequence. We will show by induction that  $S_{2k+1}\neq0$ for $k\geq 0$. In fact, for $k=1$ the equation \eqref{S2kp1} reads  
 \begin{equation*}\label{S2kp10}
-a_{1}a_n(a_{n-1}-a_{0})\,S_{3}+(a_n\omega_2-a_{1}\omega_{n+1})\,S_{1}=0,\;\, n\geq 2.
\end{equation*}
So $S_3\neq 0$, otherwise we have
\begin{equation*}\label{S2kp1011}
a_n=\frac{a_1}{\omega_2}\omega_{n+1},\,n\geq 2,
\end{equation*}
which is not solution of \eqref{S1eqan}. Suppose that $S_{2m+1}\neq  0$. Then  \eqref{S2kp1},  for $k=m+1$,  becomes
\begin{equation*}\label{S2kp11}
-a_{1}a_n(a_{n-1}-a_{0})\,S_{2m+3}+(a_n\omega_2-a_{1}\omega_{n+1})\,S_{2m+1}=0,\, n\geq 2,
\end{equation*}
leading to $S_{2m+3}\neq 0.$\\
Now, the operation $S_1\,E(n,k)-S_{2k-1}E(n,1)$ gives
\begin{equation*}
a_{n-1}\left(  \left( S_1S_{2\,k+1}-S_3S_{2\,k-1} \right) a_{ n-2} -S_1R_{2\,k+1}+R_3S_{2\,k-1} \right)=0 , \,k\geq 1,\;n\geq 2. 
\end{equation*}
Since $a_n$ is not constant, $S_1S_{2\,k+1}-S_3S_{2\,k-1}=0$  and  $-S_1R_{2\,k+1}+R_3S_{2\,k-1}=0$ for $k\geq 1$. \\
Thus, for $k\geq 1$, we obtain 
\begin{equation}\label{S2K1R2K1}
S_{2k+1}=S_1\,\left(\frac{S_3}{S_1}\right)^k\; \text{ and } R_{2k+1}=R_3\,\left(\frac{S_3}{S_1}\right)^{k-1}.
\end{equation}
 As a consequence, all equations $E(n,k)$, $k\geq 1$, will be reduced to equation $E(n,1)$. That is
 \begin{equation}\label{Skeqan3}
-S_{3}\,a_{n-2}\,a_{n-1}-S_{1}\,\omega_{n}+R_{3}\,a_{n-1}=0,\; n\geq 2.
\end{equation}
  Denote equation \eqref{S1eqan} by $E(n)$. Then   the operation $S_{3}E(n)-S_1\,E(n,1)$ gives
\begin{equation}\label{S2kpAAA}
\left( S_1R_3-S_3R_1 \right) a_{ n-1} =\frac n4\,{\frac {   \left({S_1}^{2} -2\,S_3 \right)  ( 
n+1 ) +2\,(\lambda-1) \left({S_1}^{2} -S_3 \right)  }{
   (   n+\lambda ) ( n-1+\lambda)}},\,n\geq 2.
\end{equation}
 Moreover, since $\lambda\neq1$,
  from \eqref{S2kpAAA} we must have $S_1R_3-S_3R_1\neq 0$ because we can not have simultaneously ${S_1}^{2}-2\,S_3 =0$ and ${S_1}^{2} -S_3=0$.\\
So, we can write
\begin{equation*}\label{S2kpSolu}
 a_{ n} =\frac{1}{S_1R_3-S_3R_1  }
\frac {n+1}4\,{\frac {   \left({S_1}^{2} -2\,S_3 \right)  ( 
n+2 ) +2\,(\lambda-1) \left({S_1}^{2} -S_3 \right)  }{
   (   n+1+\lambda ) ( n+\lambda)}},\,n\geq 1. 
\end{equation*}
Also notice that  $S_3\not=\frac{{S_1}^{2} }{2} $ otherwise $a_n\to0$ as $n\to\infty$, which gives a contradiction when taking the limit $n\to\infty$ in \eqref{S1eqan}.\\
So, if we set 
\begin{equation*}
\rho=\frac{{S_1}^{2} -2\,S_3 }{4(S_1R_3-S_3R_1)  }
\end{equation*}
 and 
 \begin{equation}\label{b}
 b=1+(\lambda-1)\frac{ {S_1}^{2} -S_3 }{{S_1}^{2} -2\,S_3}
 \end{equation}
then \begin{equation}\label{anlamnda}
a_{ n} =\rho
{\frac {  (n+1)  ( 
n+1 +b)  }{
   (   n+1+\lambda ) ( n+\lambda)}},\,n\geq 1. 
\end{equation}
The substitution of \eqref{anlamnda} in \eqref{S1eqan} and \eqref{Skeqan3} with the partial fraction decomposition of the resulting equations with respect to $n$  give us the following conditions
\begin{eqnarray}
&&2\,S_{{1}}{
\rho}^{2}-2\,R_{{1}}\rho+1=0,\label{first}\\
&&4\,S_{{3}}{\rho}^{2}-4\,R_{{3}}\rho +S_{{1}}=0,\\
&&\left( b-\lambda \right)  \left( - S_{{1}}\rho{\lambda}
^{2}+\,S_{{1}}\rho(b-2)\lambda +S_{{1}}\rho(b-1) +2\,R_{{1}} \right)=0,\\
&&2\,S_{{3}}{\rho}^{2}{\lambda}^{3}-4\,S_{{3}} {\rho}^{2}\left( b-1 \right) {
\lambda}^{2}+ \left( 2\,S_{{3}}{\rho}^{2} \left( {b}^{2}-3\,b+1 \right) -
4\,R_{{3}}\rho-S_{{1}} \right) \lambda\nonumber\\
&&\qquad\qquad\qquad +2\,S_{{3}}{\rho}
^{2}b \left( b-1 \right) +4\,R_{{3}}\rho b+S_{{1}}=0,\\
&&\left( \lambda -1 \right)  \left( b-\lambda +2 \right)  \left( b-
\lambda +1 \right)=0,\\
&&\left( b-\lambda +1 \right)  \left( b-\lambda \right) =0.\label{last}
\end{eqnarray}
Remark that $b\neq \lambda$. Otherwise, by \eqref{b} we have $S_3=0$ which contradicts the fact that $S_3\neq 0$. Then $b=\lambda-1$ from \eqref{last}  and  \eqref{anlamnda} simplifies to
\begin{equation}\label{blambda-1}
a_n=\rho
{\frac {  n+1  }{
     n+1+\lambda }},\,n\geq 1.
\end{equation}
As $a_n=\alpha_{n}/\alpha_{n+1}$, we find from \eqref{blambda-1} that
\begin{equation*}
\alpha_n=
\frac{\alpha_1}{\lambda +1}\frac{(\lambda+1)_n}{n!\rho^{n-1}},\,\,n\geq 2.
\end{equation*}
So, we can write
 \begin{equation*}
F( t) =
1+\frac{\alpha_1\rho}{\lambda +1} \left(  \left( 1-\frac t \rho\right) 
^{-\lambda-1}-1 \right). 
\end{equation*}
Moreover, the six equations (\ref{first})-(\ref{last}) reduce to
\begin{eqnarray}
&&2\,S_{{1}}{
\rho}^{2}-2\,R_{{1}}\rho+1=0,\label{eqx1}\\
&&4\,S_{{3}}{\rho}^{2}-4\,R_{{3}}\rho +S_{{1}}=0,\label{eqx2}\\
&&-\left( \lambda+1 \right) S_{{1}}\rho +\,R_{{1}}=0,\label{eqx3}\\
&&4\,S_{{3}} \left( \lambda+1 \right) {\rho}^{2}-4\,R_{{3}}\rho-S_{{1}}
 \left( \lambda-1 \right)
=0.\label{eqx4}
\end{eqnarray}
By \eqref{eqx1} and \eqref{eqx3} we remark that $\lambda\neq 0$.
The resolution of the system of the last four equations gives
\begin{equation*}
R_{{1}}={\frac {1+\lambda}{2\lambda\,\rho}},\;R_{{3}}={
\frac {1}{4\lambda\,{\rho}^{3}}},\;S_{{1}}={\frac {1}{2\lambda\,{\rho}^{2}}
},\;S_{{3}}={\frac {1}{8\lambda\,{\rho}^{4}}}.
\end{equation*}
Then, in virtue of \eqref{S2K1R2K1}, we get
\begin{equation*}\label{S2K1R2K112}
S_{2k+1}=\frac{2}{\lambda}\left(\frac{1	}{4\rho^2}\right)^{k+1},\; k\geq 0,
\end{equation*}
and
\begin{equation*}\label{S2K1R2K113}
 R_{2k+1}=\frac{2}{\lambda}\left(\frac{1}{2\rho}\right)^{2k+1},\; k\geq 1,
\end{equation*}
which lead to
\begin{equation*}
A(t)=\left(1-\frac{t^2}{4\rho^2}\right)^{-\frac{1}{\lambda}}\text{ and } R(t)=\rho\left(\left(1+\frac{t^2}{4\rho^2}\right)A(t)-1\right).
\end{equation*}
The choice   $\alpha_1=\frac{\lambda +1}{\rho}$ gives simply 
$$F( t) =  \left( 1-\frac t \rho\right) 
^{-\lambda-1}.$$
Then we find that
\begin{eqnarray*}\label{GFultras1}
	A(t)F(xtA(t)-R(t))&=&\left(1-\frac{t^2}{4\rho^2}\right)\left(1- \frac{xt}{\rho}+\frac{t^2}{4\rho^2}\right)^{-\lambda-1}\\
	&=&\sum_{n=0}^{+\infty}\frac{(\lambda+1)_n}{\rho^nn!} C_n^{(\lambda)}(x)t^n.  
\end{eqnarray*}

 \textbf{ii) } If $\lambda=1$,   let $a$ such that
\begin{equation}
-S_1\,a^2+R_{1}\,a=\frac{1}{2}.\label{a1} 
\end{equation}
$\bullet$ If $a_n=a$ for $n\geq 1$  then for $n=2$ equation \eqref{S1eqan} becomes
\begin{equation}
-S_1\,a_0\,a+R_{1}\,a=\frac{1}{2},\label{a01}
\end{equation}
and then the substraction Eq\eqref{a1}-Eq\eqref{a01} leads to
\begin{equation*}
-S_1\,(a-a_0)\,a=0
\end{equation*}
giving $a_0=a$. 
Notice that if $(a_n)_{n\geq 0}$ is a constant sequence, solution of \eqref{S1eqan}, then we must have  $a_n=a$ for $n\geq 0$.
So,
\begin{equation*}
F(t)=\frac{a}{a-t}.
\end{equation*}
 Furthermore, equation \eqref{Skeqanshft} reads
\begin{equation*}
-S_{2k+1}\,a^2-\frac{S_{2k-1}}{4}+R_{2k+1}\,a=0,\;k\geq 1,\label{a2}
 \end{equation*}
  which leads to
\begin{equation*}
-a^2\left(\frac{A'(t)}{A(t)}-S_1t\right)-\frac{t^2}{4}\,\frac{A'(t)}{A(t)}+a\left(\frac{R'(t)}{A(t)}-R_1t\right)=0.
\end{equation*}
In virtue of \eqref{a1} we get 
\begin{equation*}\label{eqdifS1c=0}
-a^2\frac{A'(t)}{A(t)}-\frac{t^2}{4}\,\frac{A'(t)}{A(t)}+a\frac{R'(t)}{A(t)}-\frac{t}{2}=0
\end{equation*}
and by taking into account that  $A(0)=1$ and $R(0)=0$ we obtain
\begin{equation*}
R(t)= \left( a+\frac{t^{2}}{4a}
 \right)\,A ( t ) -a.
\end{equation*}
Then we have
\begin{equation*}
A(t)F(xtA(t)-R(t))={\frac {1}{1-\frac{xt}{a}+\frac{{t}^{2}}{4a^2}}}=\sum_{n=0}^{+\infty}\frac{U_n(x)}{a^n}t^n
\end{equation*}
where $A(t)$ is arbitrary.\\
$\bullet$ If  $(a_{n})_{n\geq 0}$ is not a constant sequence then $S_3\neq 0$. Otherwise, if $S_3=0$, then for $k=1$ equation \eqref{S2kp1} reads  
 \begin{equation*}
\frac{1}{4}(a_n-a_{1})\,S_{1}=0,\;\, n\geq 2.
\end{equation*}
Thus $a_n=a_{1}=a$ for $n\geq 2$ which contradicts the hypothesis.\\
Furthermore, since $(a_{n})_{n\geq 0}$ is not a constant sequence, then there exists $n_0\geq 1$ such that $a_{n_0}\neq a_{0}$. So,   for $k=1$ and $n=n_0+1$, equation \eqref{S2kp1} becomes  
 \begin{equation*}\label{S2kp10a}
-a_{1}a_{n_0+1}(a_{n_0}-a_{0})\,S_{3}+\frac{1}{4}(a_{n_0+1}-a_{1})\,S_{1}=0,
\end{equation*}
from which we must have $a_{n_0+1}\neq a_{1}$. For $n=n_0+1$, equation \eqref{S2kp1}  leads to
\begin{equation*}\label{S2kp110}
-a_{1}a_{n_0+1}(a_{n_0}-a_{0})\,S_{2k+1}+\frac{1}{4}(a_{n_0+1}-a_{1})\,S_{2k-1}=0,\,k\geq1.
\end{equation*}
So, we can easily show by induction  that $S_{2k+1}\neq 0$ for $k\geq 1$.\\
Now, the operation $S_{2k+1}E(n)-S_1\,E(n,k)$ gives
\begin{equation}\label{S2kp333}
\left( S_{2\,k+1}R_1-S_1R_{2\,k+1} \right) a_{ n-1} +\,\frac{1}{4}(-2\,S_{2\,k+1}+S_1S_{2\,k-1})=0
, \,k\geq 1,\;n\geq 2. 
\end{equation}
Denote equation \eqref{S2kp333} by $F(n,k)$ and make the operation $F(n_0+2,k)-F(2,k)$ to get 
\begin{equation*}
\left( S_{2\,k+1}R_1-S_1R_{2\,k+1} \right)( a_{ n_0+1}-a_1)=0
, \,k\geq 1.
\end{equation*}
Then, for $k\geq 1$, $ S_{2\,k+1}R_1-S_1R_{2\,k+1}=0$ and $-2\,S_{2\,k+1}+S_1S_{2\,k-1}=0$ by using \eqref{S2kp333}. So, we can write
\begin{equation}\label{S2K1R2K123}
S_{2k+1}=S_1\left(\frac{S_1}{2}\right)^k\,\text{ and } R_{2k+1}=R_1\left(\frac{S_1}{2}\right)^k,\,k\geq 1,
\end{equation}
giving
\begin{equation*}
A(t)=\frac{2}{2-S_1t^2}\,\text{ and } R(t)=\frac{R_1t^2}{2-S_1t^2}.
\end{equation*}
Furthermore, if $\lambda=1$ and in virtue of \eqref{S2K1R2K123} the equations \eqref{S1eq} and \eqref{Skeqk2} reduce to the same equation, namely
\begin{equation*}\label{S1eqb0}
-S_1\frac{\alpha_{n-2}}{\alpha_n}+R_{1}\frac{\alpha_{n-1}}{\alpha_n}=\frac{1}{2}, n\geq 2,
\end{equation*}
which can be written as
 \begin{equation}\label{difeqb00}
-2\,S_1\,\alpha_{n-2}+\,2\,R_{1}\,\alpha_{n-1}-\,\alpha_n=0, n\geq 2.
\end{equation}
From this difference equation  and by using the generating function of the sequence $\{\alpha_n\}$, i.e. $F(t)$, we can easily obtain the  first order differential equation,
\begin{equation*}
\left( 1-2R_{{1}}t+2S_{{1}}{t}^{2} \right) {\frac {d}{dt}}F \left( t
 \right) + \left( -2R_{{1}}+4S_{{1}}t \right) F \left( t \right) +2
R_{{1}}-\alpha_{{1}}=0,
\end{equation*}
which has the following solution
\begin{equation*}
F ( t ) ={\frac { 1+\left(\alpha_1 -2R_1 \right) t
}{1-2R_1t+2S_1{t}^{2}}}.
\end{equation*}
Notice that the solution of the  difference equation \eqref{difeqb00} is given by
\begin{eqnarray*}
\alpha_n=(8S_1)^{n/2}U_n\left(\frac{R_1}{\sqrt{ 2S_1}}\right)+(\alpha_1-2R_1)(8S_1)^{(n-1)/2}U_{n-1}\left(\frac{R_1}{\sqrt{ 2S_1}}\right), \,n\geq 1,
\end{eqnarray*}
which are orthogonal polynomials in the variable $R_1$ by taking $\alpha_{{1}}$ a polynomial of degree one in $R_1$. Specifically, for $\alpha_1=2R_1$ we meet the Chebychev polynomials of the  second kind.\\
Finally, the desired GF has the form
\begin{eqnarray*}
&&A(t)F(xtA(t)-R(t))=\\
&=&{\frac { 4-2\,S_{{1}}{t}^{
2}+2(\alpha_1-2R_1)\left(2\,xt -\,R_{{1}}{t}^{2} \right) }{4-8\,R_{{1}}xt+ 4\left( 2\,{x}^{2}S_{{1}}+\,{R_{{1}}}^
{2}-\,S_{{1}} \right) {t}^{2}-4\,R_{{1}}S_{{1}}x{t}^{3}+{S_{{1}}}
^{2}{t}^{4}}}\nonumber\\
&=&1+ \sum_{n=1}^{+\infty}\alpha_n U_n(x)t^n.
\end{eqnarray*}
 
\textbf{Case 2.  $S_1= 0$.}  Then $R_1\neq 0$ and the equation \eqref{S1eqan} gives
\begin{equation}\label{S1eq0}
a_{n-1}=\frac{1}{R_{1}}\frac{n}{2(n-1+\lambda)}, n\geq 2.
\end{equation}
Substitute \eqref{omg1} and \eqref{S1eq0} in \eqref{Skeqanshft}  to get for $k\geq 1$
\begin{equation*}
\frac{N(n,k)}{16(n+\lambda-2)(n+\lambda -1)(n+\lambda)}=0,
\end{equation*}
where
\begin{eqnarray*}
N(n,k)&=&4[S_{2k+1}+R_1^{
2}S_{2k-1}-2R_1R_{2k+1}]n^3+\\
&&4(\lambda-1)[S_{2k+1}+3R_1^{2}S_{2k-1}-4R_1
R_{2k+1}]n^2-\\
&&4[\lambda S_{2k+1}- (2\lambda-1)
  (\lambda-2) R_1^{2}  S_{2k-1}+2\lambda(\lambda-2)
R_1 R_{2k+1}]n.
\end{eqnarray*}
By equating coefficients of power of $n$ in the numerator  with $0$, we obtain for $k\geq 1$
\begin{eqnarray}
&&S_{2k+1}+R_1^{
2}S_{2k-1}-2R_1R_{2k+1}=0,\label{EQSR1}\\
&&(\lambda-1)[S_{2k+1}+3R_1^{2}S_{2k-1}-4R_1
R_{2k+1}]=0,\label{EQSR2}\\
&&\lambda S_{2k+1}- (2\lambda-1)
  (\lambda-2) R_1^{2}  S_{2k-1}+2\lambda(\lambda-2)
R_1 R_{2k+1}=0.\label{EQSR3}
\end{eqnarray}

We will consider the two sub-cases, $\lambda\neq 1$ and $\lambda=1$. 
 
\textbf{i) } If $\lambda \neq 1$  then by combining the three equations \eqref{EQSR1}, \eqref{EQSR2} and \eqref{EQSR3} we can easily show by induction that $S_{2k+1}=0$ for $k\geq 0$ and $R_{2k+1}=0$ for $k\geq 1$. So, $A(t)=1$ and $R(t)=\frac{R_1}{2}t^2$.\\
Moreover,  \eqref{S1eq0} is equivalent to 
\begin{equation*}
\frac{\alpha_{n-1}}{\alpha_n}=\frac{1}{R_1}\frac{n}{2(n-1+\lambda)},\;n\geq 2,
\end{equation*}
giving for $n\geq 2$ the expression
\begin{equation*}
\alpha_n=\left\{
\begin{array}{ll}
\frac{\alpha_1}{\lambda}\frac{(\lambda)_n(2R_1)^{n-1}}{n!}&if\,\,\lambda\neq 0,\\
\alpha_1\frac{(2R_1)^{n-1}}{n}&if\,\,\lambda= 0.
\end{array}\right.
\end{equation*}
So, we can write
 \begin{equation*}
F( t) =\left\{\begin{array}{ll}
1+\frac{\alpha_1}{2\lambda R_1} \left(  ( 1-2\,R_1\,t ) 
^{-\lambda}-1 \right) & if\,\,\lambda\neq 0,\\
1-\frac{\alpha_1}{2R_1}\ln(1-2R_1 t)&if\,\,\lambda= 0.
\end{array}
\right.
\end{equation*}
For $\lambda\neq 0$, the choice $\alpha_1=2\lambda R_1$ gives $F( t) = ( 1-2\,R_1\,t ) 
^{-\lambda}$ and then
\begin{equation*}
A(t)F(xtA(t)-R(t))=( 1-2\,R_1\,xt +R_1^2t^2)^{-\lambda}=\sum_{n=0}^{+\infty}\frac{2^nR_1^n(\lambda)_n}{n!} C_n^{(\lambda)}(x)t^n.
\end{equation*}
For $\lambda= 0$, we can choose $\alpha_1=2 R_1$ to get $F( t) = 1-\ln(1-2R_1 t)$ and 
\begin{equation*}
A(t)F(xtA(t)-R(t))=1-\ln( 1-2\,R_1\,xt +R_1^2t^2)=1+\sum_{n=1}^{+\infty}\frac{2^nR_1^n}{n} T_n(x)t^n.
\end{equation*}

\textbf{ii) } If $\lambda=1$ then  $a_n=\frac{1}{2R_1}$ for $n\geq 1$. Furthermore the  three  equations \eqref{EQSR1}, \eqref{EQSR2} and \eqref{EQSR3} reduce to
\begin{equation*}
S_{2\,k+1}+R_1^{
2}S_{2\,k-1}-2R_1\,R_{2\,k+1}=0,
\end{equation*}
which leads to the differential equation
\begin{equation*}\label{eqdifc=0}
\frac{A'(t)}{A(t)}+R_1^{
2}\,t^2\frac{A'(t)}{A(t)}-2R_1\,\left(\frac{R'(t)}{A(t)}-R_1t\right)=0.
\end{equation*}
By taking into account that  $A(0)=1$ and $R(0)=0$ we obtain
\begin{equation}\label{RA}
R(t)={\frac {  \left( 1+R_1^{2}\,t^2
 \right)\,A ( t ) -1}{2R_1}}.
\end{equation}
We have the following two cases:

$\bullet$ If $a_0=\frac{1}{2R_1}$ then $$F(t)=\frac{1}{1-2R_1t}.$$
In this case $A(t)$ is arbitrary and we have
\begin{equation*}
A(t)F(xtA(t)-R(t))=\frac{1}{1-2\,txR_{{1}}+{R_{{1}}}^{2}{t}^{2} }=\sum_{n=0}^{+\infty}2^nR_1^nU_n(x)t^n
\end{equation*}

$\bullet$ If $a_0\neq\frac{1}{2R_1}$ (i.e $\alpha_1\neq 2R_1$) then  $$F ( t ) ={\frac { 1+\left(\alpha_1 -2R_1 \right) t
}{1-2R_1t}}.$$
Moreover equation \eqref{S2kp1} can be written as
\begin{equation*}
-\frac{1}{4R_1^2}\left(\frac{1}{2R_1}-a_{0}\right)\,S_{2k+1}=0, k\geq 1,
\end{equation*}
which gives $S_{2k+1}=0$ for $k\geq 1$. Then $A(t)=1$ and $R(t)=\frac{R_1}{2}t^2$ by \eqref{RA}. In this case the GF takes the form
\begin{eqnarray*}
A(t)F(xtA(t)-R(t))&=&{\frac { 1+\left(\alpha_1 -2R_1 \right) (xt-\frac{R_1}{2}t^2)
}{1-2R_1xt+{R_{{1}}}^{2}{t}^{2}}}\nonumber\\
&=&1+\alpha_1\sum_{n=1}^{+\infty}2^{n-1}R_1^{n-1}U_n(x)t^n.
\end{eqnarray*}
\end{proof}

\subsection{Hermite polynomials}
The monic Hermite polynomials ${H}_{n}(x)$ are defined by
\begin{equation*}
H_n(x)=x^n\, {}_2F_0\left(
	\begin{array}{cr}
		-\frac{n}{2}\,\,\,\,-\frac{n-1}{2}&\\
		&;-\frac{1}{x^2}\\
		-&
	\end{array}\right), \,n\geq 0.
\end{equation*}
In this case we need to recall the  associated monic Hermite polynomials 
\begin{equation}\label{AssociateH}
		H_n(x,c)=\sum_{k=0}^{[\frac{n}{2}]}(-2)^{-k}\frac{(c)_{k}(n-k)!}{k!(n-2k)!}{H}_{n-2k}\left(x\right), \,n\geq 1, 
	\end{equation}
	which satisfy the three term difference equation
	\begin{eqnarray}\label{rec_as_herm}
		&&xH_{n}(x,c)=H_{n+1}(x,c)+\frac{n+c}{2}H_{n-1}(x,c),\;\;n\geq 0\\
		&& H_{-1}(x,c)=0,\;\;H_{0}(x,c)=1.\nonumber
	\end{eqnarray}
	and are orthogonal with respect to a
	positive measure when $c>-1$. Some properties of the polynomials $H_{n}(x,c)$ were given in  \cite{AskeyWimp1984}.
\begin{theorem}\label{th2}
	The only generating functions of the form \eqref{gf0} of  the monic Hermite polynomials $H_n(x)$ are
	\begin{eqnarray}	   
		&&\bullet \sum_{n=0}^{+\infty}\frac{1}{n!} H_n(x)t^n=e^{tx-\frac{t^2}{4}}, \mbox{ with } A(t)=1,\; R(t)=\frac{t^2}{4} \mbox{ and } F(t)=e^x.\label{H1}\\
		&&\bullet\sum_{n=0}^{+\infty}\frac{2^n}{n!}[H_n(\rho)+(\tfrac{\alpha_1}{2}-\rho)H_{n-1}(\rho,1)] H_n(x)t^n\label{H2}\\
		&&\qquad=\tfrac{1}{\sqrt{1-t^2}}\left(1+ (\tfrac{\alpha_1}{2}-\rho)\sqrt{\pi}{\rm e}^{-{{\rho}^{2}}}   \left[{\rm erfi} \left(\rho\right)-{\rm erfi} \left({\tfrac {\rho-tx}{\sqrt {1-{t}^{2}}}}\right) \right]\right)e^{\tfrac{2x\rho t-(x^2+\rho^2)t^2}{1-t^2}}\;\;\;\nonumber\\
		&& \mbox{ with }A(t)=\tfrac{1}{\sqrt{1-\,t^2}},\;  R(t)=\rho\left(\tfrac{1}{\sqrt{1-\,t^2}}-1\right)\nonumber\\
		&& \mbox{ and } F ( t)= \left(1+  (\tfrac{\alpha_1}{2}-\rho)\sqrt {\pi}{\rm e}^{-\,
			{\rho}^{2}}\left[{\rm erfi} \left(\rho\right)-{\rm erfi} \left(\rho- t\right) \right]\right)
		{ e}^{t\left( 2\rho-t \right)} ,\nonumber
	\end{eqnarray}
	where  $\rho$ is a complex number.
	 \end{theorem}	
\begin{remark} From \eqref{H2} we obtain the well known Mehler's formula
\begin{eqnarray}
&&\sum_{n=0}^{+\infty}\frac{2^n}{n!}H_n\left(\rho\right) H_n(x)t^n=\frac{1}{\sqrt{1-t^2}}e^{\frac{2x\rho t-(x^2+\rho^2)t^2}{1-t^2}}
\end{eqnarray}
 and  the following bilateral GF 
	\begin{eqnarray}
&&\sum_{n=1}^{+\infty}\frac{2^n}{n!}H_{n-1}(\rho,1) H_n(x)t^n=\frac{\sqrt{\pi}}{\sqrt{1-t^2}}\left({\rm erfi} \left(\rho\right)-{\rm erfi} \left({\tfrac {\rho-tx}{\sqrt {1-{t}^{2}}}}\right)\right)e^{-\tfrac{(\rho-tx)^2 }{1-t^2}}.\;\;\;\label{H44}\nonumber\\
\end{eqnarray}
	 
\end{remark}
	\begin{proof}[Proof of Theorem~\ref{th2}]
The  Jacobi-Szeg\"o parameters of the monic Hermite polynomials are given by
  $\beta_n=0$ (symmetric) and   
\begin{equation}\label{omg1_1}
\omega_n=\frac {n }{ 2 },\;n\geq1.
\end{equation}
Then by  \eqref{eqiis} we obtain
\begin{equation}\label{solA1_1}
A^1_n=\frac{n}{2}, \;n\geq 2,
\end{equation}
and the equation \eqref{A1n11} can be written as
\begin{equation}\label{S1eqAA}
-S_1\frac{\alpha_{n-2}}{\alpha_n}+R_{1}\frac{\alpha_{n-1}}{\alpha_n}=\frac{n}{2}, n\geq 2.
\end{equation}
By \eqref{solA1_1} and \eqref{eqivs} we get $A_n^3=0$ for $n\geq 4$ and according to Corollary \ref{corl1} we have 
$$A_n^{2k+1}=0 \text{ for } k\geq 1$$
which leads to
\begin{equation}\label{Skeqk2_2}
-S_{2k+1}\,\frac{\alpha_{n-2}}{\alpha_{n}}-S_{2k-1}\,\omega_{n}\,+R_{2k+1}\,\frac{\alpha_{n-1}}{\alpha_{n}}=0,k\geq 1,\, n\geq 2
\end{equation}
after multiplying \eqref{Aknnn} by $\frac{\alpha_n}{\alpha_{n-2k}}$ and using the shifting $n\to n+2k$.\\ 
The equations \eqref{S1eqAA} and \eqref{Skeqk2_2} are equivalent, respectively, to
\begin{equation}\label{S1eqan1}
-S_1\,a_{n-2}\,a_{n-1}+R_{1}\,a_{n-1}=\frac{n}{2}, n\geq 2,
\end{equation}
and
   \begin{equation}\label{Skeqanshft1}
-S_{2k+1}\,a_{n-2}\,a_{n-1}-S_{2k-1}\,\omega_{n}+R_{2k+1}\,a_{n-1}=0,\;k\geq 1,\, n\geq 2.
\end{equation}
Observe that from \eqref{S1eqan1} we must have $(S_1,R_1)\neq (0,0)$. \\

  \textbf{Case 1.  $S_1\neq 0$.}  By analogy with the case of ultraspherical polynomials we show that, for $k\geq 1$, 
\begin{equation}\label{S2K1R2K1_1}
S_{2k+1}=S_1\,\left(\frac{S_3}{S_1}\right)^k\; \text{ and } R_{2k+1}=R_3\,\left(\frac{S_3}{S_1}\right)^{k-1},
\end{equation}
 and all of the equations \eqref{Skeqanshft1}, for $k\geq 1$, will be reduced to 
 \begin{equation}\label{Skeqanshft1k3}
-S_{3}\,a_{n-2}\,a_{n-1}-S_{1}\,\omega_{n}+R_{3}\,a_{n-1}=0,\; n\geq 2.
\end{equation}
  Now the operation $S_{3}$Eq\eqref{S1eqan1}-$S_1$\,Eq\eqref{Skeqanshft1k3} gives
\begin{equation}\label{S2kpBBB_1}
\left( S_1R_3-S_3R_1 \right) a_{ n-1} =\frac {  -S_3+{S_1}^{2} }{2}\,n ,\,n\geq 2.
\end{equation}
Then  we must have $S_1R_3-S_3R_1=0$ and $-S_3+{S_1}^{2}=0$, otherwise $a_n$ given by \eqref{S2kpBBB_1} cannot be a solution of the system of equations \eqref{S1eqan1} and \eqref{Skeqanshft1k3}. So, in virtue of \eqref{S2K1R2K1_1} we can write 
\begin{equation}\label{S2K1R2K12_1}
S_{2k+1}=S_1^{k+1}\; \text{ and } R_{2k+1}=R_1\,S_1^{k} \;, k\geq 1,
\end{equation}
which leads to the first order differential equations
\begin{equation}\label{ATRT_1}
\frac{A'(t)}{A(t)}=\frac{S_1\,t}{1-S_1\,t^2}\;\text{ and } \;\frac{R'(t)}{A(t)}=\frac{R_1\,t}{1-S_1\,t^2}.
\end{equation}
Now by solving \eqref{ATRT_1} we obtain
\begin{equation*}
A(t)=\frac{1}{\sqrt{1-S_1\,t^2}}\;\text{ and } R(t)=\frac{R_1}{S_1}\left(\frac{1}{\sqrt{1-S_1\,t^2}}-1\right).
\end{equation*}
Furthermore, according to \eqref{S2K1R2K12_1}, the equations \eqref{S1eqan1} and \eqref{Skeqanshft1k3} reduce to the same equation, namely
 \begin{equation}\label{difeqb0}
-2\,S_1\,\alpha_{n-2}+\,2\,R_{1}\,\alpha_{n-1}-n\,\alpha_n=0,\; n\geq 2.
\end{equation}
From this difference equation  and by using the generating function of the sequence $\{\alpha_n\}$, i.e. $F(t)$, we can easily obtain the following  first order differential equation   
\begin{equation*}\label{DEbeq0}
{\frac {d}{dt}}F ( t) + 2\,(-  R_1+  S_1 t ) F ( t ) + 2\,R_1 -\alpha_1=0
\end{equation*}
which has the following solution
\begin{eqnarray*}
F ( t)& =& \left(1+  \left(\alpha_1 -2 R_1
 \right) \frac{\sqrt {\pi}}{2\sqrt {S_1}}{\rm e}^{-\,\frac {
{R_1}^{2}}{S_1}}\left[{\rm erfi} \left(\frac{R_1}{\sqrt {\,S_1}}\right)\right.\right.\\
&&\qquad\qquad\qquad\qquad\quad\qquad\qquad\qquad\left.\left.-{\rm erfi} \left(\frac{ R_1-S_1 t}{\sqrt {{ S_1}}}\right) \right]\right)
{ e}^{t\left( 2R_1-S_1t \right)} .
\end{eqnarray*}
If we put $\alpha_n=\frac{(4S_1)^{\frac{n	}{2}}}{n!}\gamma_n$ then after shifting $n\to n+1$, \eqref{difeqb0} becomes 
\begin{equation}\label{difeqb001}
	\rho\gamma_{n}=\gamma_{n+1}+\frac{n}{2}\,\gamma_{n-1},\; n\geq 1.
\end{equation}
with $\rho=\frac{R_1}{\sqrt{S_1}}$. Remark that $H_n(\rho)$ and $H_{n-1}(\rho,1)$ are  two independent  solutions of \eqref{difeqb001}. So the solution of  \eqref{difeqb0} has the form
\begin{eqnarray*}
\alpha_n =\frac{(4S_1)^{\frac n2}}{n!}\left( H_n\left(\rho\right)
+\frac{\alpha_1-2R_1}{2\sqrt{S_1}}H_{n-1}(\rho,1)\right) , \,n\geq 1,
\end{eqnarray*}
which are orthogonal polynomials in the variable $R_1$ by taking $\alpha_{{1}}$ a polynomial of degree one in $R_1$. Clearly for $\alpha_1=2R_1$ the $\alpha_{n}$ are Hermite polynomials.\\
The GF $A(t)F(xtA(t)-R(t))$ takes the form
\begin{eqnarray*}
&&\tfrac{1}{\sqrt{1-S_1\,t^2}}\left(1+  \tfrac{\left(\alpha_1-2 R_1
	\right)\sqrt {\pi}}{2\sqrt {S_1}}{\rm e}^{-\,\frac{{R_1}^{2}}{S_1}}
\left[{\rm erfi}\left(\tfrac{R_1}{\sqrt{S_1}}\right)
 \right.\right.
 \\&&\qquad\qquad\qquad\left.\left.- {\rm erfi} \left(\tfrac{R_1-S_1xt}{\sqrt{S_1}\sqrt{1-S_1t^2}}
\right) \right]\right)
{ e}^{-\frac{S_1x^2t^2-2\,R_1xt+{R_1}^2t^2}{1-S_1t^2}}=\sum_{n=0}^{+\infty}\alpha_n H_n(x)t^n.
\end{eqnarray*}

\textbf{Case 2.  $S_1= 0$.} Then  $R_1\neq 0$ and the equation \eqref{S1eqan1} gives
\begin{equation}\label{S1eq0_3}
a_{n-1}=\frac{n}{2R_{1}}, n\geq 2.
\end{equation}
We substitute \eqref{S1eq0_3} in \eqref{Skeqanshft1} and we show by induction that $S_{2k+1}=0$ for $k\geq 0$ and $R_{2k+1}=0$ for $k\geq 1$. So, $A(t)=1$ and $R(t)=\frac{R_1}{2}t^2$.
Moreover,  \eqref{S1eq0_3} is equivalent to 
\begin{equation*}
\frac{\alpha_{n-1}}{\alpha_n}=\frac{n}{2R_{1}},\;n\geq 2
\end{equation*}
giving, for $n\geq 2$,
\begin{equation*}
\alpha_n=\alpha_1\frac{(2\,R_1)^{n-1}}{n!}
\end{equation*}
and
 \begin{equation*}
F( t) =1+\frac{\alpha_1}{2\,R_1}(e^{2\,R_1\,t}-1).
\end{equation*}
The choice $\alpha_1=2 R_1$ gives  $F( t) = e^{2\,R_1\,t}$ and then
\begin{equation*}
A(t)F(xtA(t)-R(t))=e^{2\,R_1\,xt-R_1^2t^2}=\sum_{n=0}^{+\infty}\frac{2^nR_1^n}{n!} H_n(x)t^n.
\end{equation*}
\end{proof}



\subsection{Jacobi polynomials}
The monic Jacobi polynomials $\tilde{P}_{n}^{(\alpha,\beta)}(x)$ are defined by
\begin{equation*}
\tilde{P}_{n}^{(\alpha,\beta)}(x)=\frac{2^n(\alpha+1)_n}{n!(n+\alpha+\beta+1)_n}{}_2F_1\left(
\begin{array}{cr}
	-n\,\,\,\,n+\alpha+\beta+1&\\
	&;\frac{1-x}{2} \\
	\alpha+1&
\end{array}\right).
\end{equation*}
\begin{theorem}
	The only generating functions of the form \eqref{gf0} of  the monic Jacobi polynomials $\tilde{P}_{n}^{(\alpha,\beta)}(x)$, ($\alpha\neq \beta$) are
	\begin{eqnarray}
		&&\bullet\sum_{n\geq0}\tfrac{\left( \alpha+\tfrac{3}{2}\right)_{n}}{n!}\tilde{P}_{n}^{(\alpha,\alpha+1)}(x)t^{n}=\left(1-\tfrac{t}{2}\right)\left(1- xt+\tfrac{t^2}{4}\right)^{-(\alpha+\tfrac{3}{2})}, \label{J1}\\
		&& \mbox{ with } 	A(t)=(1-\tfrac{ t}{2})^{\tfrac{-2}{2\alpha+1}},\, R(t)=\left(1+\tfrac{t^2	}{4} \right)A(t)-1 \text{ and } F(t)=\left(1- t \right)^{-(\alpha+\tfrac{3}{2})}.  \nonumber\\ 
		&&\bullet\sum_{n\geq0}\tfrac{\left(\tfrac{3}{2}\right)_{n}}{(n+1)!}\tilde{P}_{n}^{(0,1)}(x)t^{n}  =\tfrac {1}{ \left( 1-t\right) ^{2}} {}_2F_1\left(
		\begin{array}{cr}
			\tfrac{3}{2}\,\,\,\,1&\\
			&;\tfrac{(x+1)t}{\left( 1+\tfrac{t}{2}\right)^{2}}\\
			2&
		\end{array}\right),\label{J2}\\
		&&\mbox{with }		A(t)=\tfrac {1}{ \left( 1-t\right) ^{2}},\, R(t)=\left(1+\tfrac{t^2}{4} \right)A(t)-1, \mbox{ and }F(t)={}_2F_1\left(
		\begin{array}{cr}
			\tfrac{3}{2}\,\,\,\,1&\\
			&;t\\
			2&
		\end{array}\right).\nonumber\\
		&&\bullet\sum_{n\geq 0}\tfrac{\left(\tfrac{3}{2}\right)_n}{\left(\tfrac{3}{2}+\tfrac{2\alpha-1}{2R_0}\right)_n}\tilde{P}_{n}^{(\alpha,1-\alpha)}(x)\left(\tfrac{2t}{R_0}\right)^n={\tfrac {1}{ \left( 1-t\right) ^{2}} }\,{}_2F_1\left(
		\begin{array}{cr}
			\tfrac{3}{2}\,\,\,\,1&\\
			&;{\tfrac {2t
					\left( x-R_{{0}} \right) }{ R_{{0}}\left(1-t\right) ^{2}}}
			\\
			\tfrac{3}{2}+\tfrac{2\alpha-1}{2R_0}&
		\end{array}\right)\nonumber\\
		\;\;\;\label{J3}\\
			&&\mbox{with }A(t)=\tfrac {1}{ \left( 1-t\right) ^{2}},\, R(t)=\tfrac{R_0}{2}\left(1+t^2 \right)A(t)-\tfrac{R_0}{2}, \nonumber\\
			&&\qquad\qquad\qquad\qquad\mbox{ and }F(t)={}_2F_1\left(
		\begin{array}{cr}
			\tfrac{3}{2}\,\,\,\,1&\\
			&;\tfrac{2}{R_0} t\\
			\tfrac{3}{2}+\tfrac{2\alpha-1}{2R_0}&
		\end{array}\right).\nonumber
	\end{eqnarray}
	
\end{theorem}
\begin{proof}
The  Jacobi-Szeg\"o parameters of Jacobi polynomials are given by
\begin{eqnarray}\label{eqj1}
	\beta_n =\frac{\beta^2-\alpha^2}{4(n+\left(\alpha+\beta\right)/2)(n+\left(\alpha+\beta\right)/2+1)}
\end{eqnarray}
and
\begin{eqnarray}\label{eqj2}
	\omega_n =\frac{n\left(n+\alpha\right)\left(n+\beta\right)\left(n+\alpha+\beta\right)}{4(n+\left(\alpha+\beta-1\right)/2)(n+\left(\alpha+\beta+1\right)/2)(n+\left(\alpha+\beta\right)/2)^2}.
\end{eqnarray}
From \eqref{eqi} we have 
\begin{equation}\label{A0jacob}
	A^{0}_{n+1}=\frac{1}{n+1}\sum_{k=0}^{n}\beta_k=\frac{\beta-\alpha}{2\,n+\alpha+\beta+2}
\end{equation}
where according to \eqref{eqj1} and the expression of $A^{0}_{n+1}$ we get for $n\geq 0$
\begin{equation*}
	-S_0\frac{\alpha_{n}}{\alpha_{n+1}}+R_0=\frac{\beta-\alpha}{2\,n+\alpha+\beta+2}.
\end{equation*}
Note that, for $\alpha=\beta$ we have the symmetric case (ultraspherical polynomials).\\
Assume $\beta\neq\alpha$. So, $S_0\neq 0$ and 
\begin{equation*}
	a_n=\frac{\alpha_{n}}{\alpha_{n+1}}={\frac {R_{{0}} \left( 2\,n+\alpha+\beta+2 \right) +\alpha-\beta}{S_{{0
			}} \left( 2\,n+\alpha+\beta+2 \right) }}
\end{equation*}
leading to
\begin{eqnarray*}
	\alpha_{n}=&&{\frac {S_{{0}} \left( 2\,n+\alpha+\beta \right) }{R_{{0}} \left( 2\,n
			+\alpha+\beta \right) +\alpha-\beta}}\alpha_{n-1}\nonumber\\
		=&&\left\{\begin{array}{ll}
		\left( {\frac {S_{{0}}}{R_{{0}}}} \right) ^{n} \frac{\left( 1+\frac{\alpha+\beta}{2}\right)_{n}}{\left( 1+\frac{\alpha+\beta}{2}+{\frac {\alpha-\beta}{2R_{{0}}}}\right)_{n}},& \text{ if } R_{0}\neq0\\
		\left( {\frac {2S_{{0}}}{\alpha-\beta}}\right) ^{n} \left( 1+\frac{\alpha+\beta}{2}\right)_{n},& \text{ if } R_{0}=0.
		\end{array}
		\right.
\end{eqnarray*}
By using \eqref{A0jacob}, \eqref{eqii}, \eqref{eqiii}, \eqref{eqiv} and in virtue of \eqref{eqj1} and \eqref{eqj2} we can show that $A_n^2=0$ for $n\geq 3$ and $A_n^3=0$ for $n\geq 4$. Then according to Corollary \ref{corl1} we have $$A_n^{k}=0 \text{ for }k\geq 2 \text{ and }n\geq k+1$$ wich  are equivalent to
\begin{equation}\label{eqj3}
-S_{k}a_{n-k-1}-\beta_{n-k}S_{k-1}-\frac{\omega_{n-k+1}}{a_{n-k}}S_{k-2}+R_{k}=0.
\end{equation}
Equation $A_{n}^{2}=0$ excludes the case $R_{0}=0$. So, equations \eqref{eqj3} become
\begin{eqnarray} \label{eqj4}
&&-{\frac{S_{{k}}R_{{0}}}{S_{{0}}}}-{\frac {S_{{k-2}}S_{{0}}}{4R_{{0}}}}+R_{{k}}\nonumber\\ 
&&-\frac{S_{{k-2}}S_{{0}}\left( \alpha-\beta\right)\left(R_{{0}}^2-1\right)\left(R_{{0}}^{2}(\alpha+\beta)^{2}-(\alpha-\beta)^2\right)}{4R_{{0}} \left( R_{{0}}^2-(\alpha-\beta)^2\right)  \left( R_{{0}} \left( 2n-2k+\alpha+\beta+2\right)+\alpha-\beta \right)}\nonumber\\
&&+\frac{S_{{k-2}}S_{{0}}((\alpha+1)^2-\beta^2)((\alpha-1)^2-\beta^2)}{8\, \left( 2n-2k+3+\alpha+\beta \right)  \left(R_{{0}}+\beta-\alpha\right) }\nonumber\\
&&+\frac{\left( \alpha^{2}-\beta^{2}\right)\left( S_{{k-2}}S_{{0}
	} \left( \alpha+\beta \right) -2\,S_{{k-1}} \right) 
}{4(2n-2k+\alpha+\beta+2)}\nonumber\\
&&+\frac{\left( \alpha-\beta\right) \left( S_{{k-1}}S_{{0}} \left( \alpha+
	\beta \right) -2\,S_{{k}} \right)
}{2S_{{0}}\left(2n-2k+\alpha+\beta \right)}\nonumber\\
&&-\frac{S_{{k-2}}S_{{0}}((\alpha+1)^2-\beta^2)((\alpha-1)^2-\beta^2)}{8\, \left( 2n-2k+1+\alpha+\beta \right)  \left( R_{{0}}-\beta+\alpha\right) }=0.
\end{eqnarray}

If $\alpha+\beta=0$ then from \eqref{eqj4} we get $S_0=0$ which contradicts $S_0\neq0$. 

If $\alpha+\beta\neq0$ then we have from \eqref{eqj4}
 $$\beta\in \left\lbrace \alpha+1,\alpha-1,-\alpha+1,-\alpha-1 \right\rbrace, $$
\begin{equation*}
	S_k=S_0 \rho^k,\,\,(k\geq0),\,\text{ with } \rho=\frac{(\alpha+\beta)S_0}{2},
\end{equation*}
and
\begin{equation*}
	R_{{k}}={\frac{S_{{k}}R_{{0}}}{S_{{0}}}}+{\frac {S_{{k-2}}S_{{0}}}{4R_{{0}}}}.
\end{equation*}
Equation \eqref{eqii} gives $R_1=R_0S_1/S_0+S_0/2R_0$, so
\begin{equation*}
	A(t)=(1-\rho t)^{\frac{-2}{\alpha+\beta}}\,\text{ and } R(t)=\left(\frac{R_0}{S_0}+\frac{S_0}{4R_0}t^2 \right)A(t)-\frac{R_0}{S_0}. 
\end{equation*}
We have also from \eqref{eqii} that 
$$\alpha+\beta+2=0 \text{ or } 1+ \left( -\beta-\alpha \right) {R_{{0}}}^{2}+ \left( \alpha-1+\beta\right)  \left( -\beta+\alpha \right) R_{{0}}=0.$$

-- If $\beta=\alpha+1$ then from \eqref{eqj4} $R_0=\frac{\pm 1}{2\alpha+1}$ with $\alpha\neq -1/2$. The fact that $\alpha=-1/2$, contradicts $S_0\neq 0$ by $A_{n}^{2}=0$. For $R_0=\frac{1}{2\alpha+1}$ all equations of Proposition~\ref{prop3} are verified and 
\begin{eqnarray*}
	\alpha_{n}=&&\left( {\frac {S_{{0}}}{R_{{0}}}} \right) ^{n} \frac{\left( 1+\frac{\alpha+\beta}{2}\right)_{n}}{\left( 1+\frac{\alpha+\beta}{2}+{\frac {\alpha-\beta}{2R_{{0}}}}\right)_{n}}\nonumber\\
	=&&\left( 2\rho\right) ^{n} \frac{\left( \alpha+\frac{3}{2}\right)_{n}}{n!}.
\end{eqnarray*}
\begin{equation*}
	F(t)=\left(1-2\rho t \right)^{-(\alpha+\frac{3}{2})}, 	A(t)=(1-\rho t)^{\frac{-2}{2\alpha+1}}\,\text{ and } R(t)=\frac{1}{2\rho}\left(1+(\rho t)^2 \right)A(t)-\frac{1}{2\rho}.  
\end{equation*}
After some calculations, and without lose of generality by taking $\rho=1/2$, we find that
\begin{eqnarray}\label{GF1}
	A(t)F(xtA(t)-R(t))&=&\left(1-\frac{t}{2}\right)\left(1- xt+\frac{t^2}{4}\right)^{-(\alpha+\frac{3}{2})}\nonumber\\
	&=&\sum_{n\geq0}\frac{\left( \alpha+\frac{3}{2}\right)_{n}}{n!}\tilde{P}_{n}^{(\alpha,\alpha+1)}(x)t^{n}.  
\end{eqnarray}
For $R_0=\frac{-1}{2\alpha+1}$ we have just the case $\alpha=0$. So, $R_0=-1$ and by taking $\rho=-1/2$, the GF will have the form
 \begin{eqnarray}\label{GF2}
 	A(t)F(xtA(t)-R(t))&=&2\frac{\left(1+\frac{t}{2}\right)\left(1-xt+\frac{t^2}{4}\right)^{-\frac{1}{2}}-1}{t(x+1)}\nonumber\\
 	&=&\sum_{n\geq0}\frac{\left(\frac{3}{2}\right)_{n}}{(n+1)!}\tilde{P}_{n}^{(0,1)}(x)t^{n}  
 \end{eqnarray}
with
\begin{equation*}
	F(t)=2\frac{\left( 1-t\right) ^{-\frac{1}{2}}-1}{t}, 	A(t)=\left( 1+\frac{t}{2}\right)^{-2}\,\text{ and } R(t)=\left(1+\frac{t^2}{4} \right)A(t)-1.  
\end{equation*}
Note that in this case
\begin{eqnarray*}
	&& F(t)={}_2F_1\left(
	\begin{array}{cr}
		\frac{3}{2}\,\,\,\,1&\\
		&;t\\
		2&
	\end{array}\right)
\end{eqnarray*}
and the GF \eqref{GF2} has also the form
\begin{eqnarray*}
	&& \left( 1+\frac{t}{2}\right)^{-2} {}_2F_1\left(
	\begin{array}{cr}
		\frac{3}{2}\,\,\,\,1&\\
		&;\frac{(x+1)t}{\left( 1+\frac{t}{2}\right)^{2}}\\
		2&
	\end{array}\right).
\end{eqnarray*}

-- If $\beta=\alpha-1$ then due to the relation $\tilde{P}_{n}^{(\alpha,\beta)}(-x)=(-1)^{n}\tilde{P}_{n}^{(\beta,\alpha)}(x)$, we can write $\tilde{P}_{n}^{(\beta+1,\beta)}(x)=(-1)^{n}\tilde{P}_{n}^{(\beta,\beta+1)}(-x)$ and use the two previous cases.

-- If $\beta=-\alpha+1$ then from \eqref{eqii} $R_{0}^{2}=1$ and all equations of proposition \eqref{prop3} are verified. We have 
\begin{eqnarray*}
	&&A(t)=(1-\rho t)^{-2},\, R(t)=\frac{R_0}{2\rho}\left(1+(\rho t)^2 \right)A(t)-\frac{R_0}{2\rho},\\
&& F(t)={}_2F_1\left(
\begin{array}{cr}
	\frac{3}{2}\,\,\,\,1&\\
	&;\frac{2\rho}{R_0} t\\
	\frac{3}{2}+\frac{\alpha-\beta}{2R_0}&
\end{array}\right)
\end{eqnarray*}
and
\begin{eqnarray*}
A(t)F(xtA(t)-R(t))&=&{\frac {1}{ \left( 1-\rho\,t\right) ^{2}} }{}_2F_1\left(
\begin{array}{cr}
	\frac{3}{2}\,\,\,\,1&\\
	&;{\frac {2\rho\,t
 \left( x-R_{{0}} \right) }{ R_{{0}}\left(1- \rho\,t\right) ^{2}}}
\\
	\frac{3}{2}+\frac{2\alpha-1}{2R_0}&
\end{array}\right)\\
&=& \sum_{n\geq 0}\frac{\left(\frac{3}{2}\right)_n\left(\frac{2\rho}{R_0}\right)^n}{\left(\frac{3}{2}+\frac{2\alpha-1}{2R_0}\right)_n}\tilde{P}_{n}^{(\alpha,1-\alpha)}(x)t^{n}.
\end{eqnarray*}
-- If $\beta=-\alpha-1$ then from \eqref{eqii} $R_{0}^{2}=1$. $R_0=1$ gives $(\alpha,\beta)=(0,-1)$ and $R_0=-1$ gives $(\alpha,\beta)=(-1,0)$. They are subclasses of \eqref{GF1} by taking $\alpha=-1$ and the use of the relation $\tilde{P}_{n}^{(0,-1)}(-x)=(-1)^{n}\tilde{P}_{n}^{(-1,0)}(x)$,
\begin{equation}\label{GF3}
	\left(1-\frac{t}{2}\right)\left(1- xt+\frac{t^2}{4}\right)^{-\frac{1}{2}}=\sum_{n\geq0}\frac{\left( \frac{1}{2}\right)_{n}}{n!}\tilde{P}_{n}^{(-1,0)}(x)t^{n}.  
\end{equation}
\end{proof}
\subsection{Laguerre polynomials}
The monic Laguerre polynomials ${L}_{n}^{(\alpha)}(x)$ are defined by 
\begin{equation*}
{L}_{n}^{(\alpha)}(x)=(-1)^n(\alpha+1)_n\,{}_1F_1\left(
	\begin{array}{cr}
		-n&\\
		&;x\\
		\alpha+1&
	\end{array}\right)
\end{equation*}
\begin{theorem}\label{thLaguerre}
	The only generating functions of the form \eqref{gf0} of the monic Laguerre polynomials $L_n^{\alpha}(x)$ are
	\begin{eqnarray}
		&&\bullet \sum_{n\geq 0}L_{n}^{(\alpha)}(x)\frac{\left(- t \right)^n }{n!}=	\left(1- t \right)^{-\alpha-1}e^{\frac{xt}{t-1}}.\label{L2}\\
		&& \mbox{ with } 	A(t)=\frac{1}{1-t},\,\, R(t)=-\alpha\ln\left(1-t \right) \mbox{ and } F(t)=e^{-t}.\nonumber\\
	&&\bullet \sum_{n\geq 0}\frac{\left(- 1\right)^n}{(\alpha+1)_{n}}L_{n}^{(\alpha)}(x)t^n=	\frac{1}{1- t}\,{}_1F_1\left(
		\begin{array}{cr}
			1&\\
			&;\frac{xt}{ t-1}\\
			\alpha+1&
		\end{array}\right),	\label{L1}\\
			&&\mbox{with }A(t)=\frac{1}{1-t},\,\, R(t)=0\mbox{ and } F(t)={}_1F_1\left(
			\begin{array}{cr}
				1&\\
				&;-t\\
				\alpha+1&
			\end{array}\right).\nonumber
	\end{eqnarray}
	
\end{theorem}
\begin{proof}
The Jacobi-Szeg\"o parameters  are given by
\begin{equation}\label{LjacSz}
	\beta_{n}=2n+\alpha+1 \text{ and } \omega_{n}=n(n+\alpha).
\end{equation}
From \eqref{eqi} we have for $n\geq 0$
\begin{equation}\label{A0Lag}
	A^{0}_{n+1}=-S_0\frac{\alpha_{n}}{\alpha_{n+1}}+R_0=\frac{1}{n+1}\sum_{k=0}^{n}\beta_k=n+\alpha+1.
\end{equation}
Remark that $S_{0}\neq 0$, so
\begin{eqnarray*}
	\alpha_{n}=&& {\frac {-S_{0}}{n+1+\alpha-R_{0}}}\alpha_{n-1}\nonumber\\
	=&& \frac{\left( -S_{0}\right) ^{n}}{\left( 1+\alpha-R_{0}\right) _{n}}.
\end{eqnarray*}
By using \eqref{A0Lag}, \eqref{eqii}, \eqref{eqiii}, \eqref{eqiv} and in virtue of \eqref{LjacSz} we can show that $A_n^2=0$ for $n\geq 3$ and $A_n^3=0$ for $n\geq 4$. Then according to Corollary \ref{corl1} we have $$A_n^{k}=0 \text{  for }k\geq 2\text{  and } n\geq k+1$$
or equivalently
\begin{eqnarray}\label{eqL1}
	&&\left( -2\,S_{{k-1}}+S_{{k-2}}S_{{0}}+{\frac {S_{{k}}}{S_{{0}}}}
	\right) n+{\frac {S_{{k}} \left( -k-R_{{0}}+\alpha \right) }{S_{{0}}}}\nonumber\\
&&-\left( k-1-R_{{0}} \right) S_{{k-2}}S_{{0}}- \left( -2\,k+1+\alpha
	\right) S_{{k-1}}+R_{{k}}\\
	&&-{\frac {R_{{0}} \left( -R_{{0}}+\alpha
			\right) S_{{k-2}}S_{{0}}}{n-k+1-R_{{0}}+\alpha}}=0.\nonumber
\end{eqnarray}
We see that $R_0=0$ or $R_0=\alpha$. 

-- If  $R_0=\alpha$ then $\alpha_n=(-S_0)^n/n!$ and $F(t)=\exp\left(-S_{0}t \right) $. From the relation of $\omega_{n}$ we get $S_1=S_{0}^2$ and $R_1=0$. Equation \eqref{eqL1} gives $S_k=S_{0}^{k+1}$ and $R_{k}=0$ for $k\geq 2$. So,
\begin{eqnarray*}
	A(t)=\frac{1}{1-S_0 t},\,\, R(t)=-\frac{\alpha}{S_0}\ln\left(1-S_0 t \right)
\end{eqnarray*}
and 
\begin{eqnarray*}
	A(t)F(xtA(t)-R(t))=\left(1-S_0 t \right)^{-\alpha-1}e^{\frac{-S_0 xt}{1-S_0 t}}=\sum_{n\geq 0}L_{n}^{(\alpha)}(x)\frac{\left(-S_0 t \right)^n }{n!}.
\end{eqnarray*}

-- If  $R_0=0$ then  $\alpha_n=(-S_0)^n/(\alpha+1)_{n}$, $S_k=S_{0}^{k+1}$ and $R_{k}=0$ for $k\geq 0$, giving $A(t)=1/(1-S_0 t)$, $R(t)=0$ and
\begin{eqnarray*}
	A(t)F(xtA(t)-R(t))&=&\frac{1}{1-S_0 t}{}_1F_1\left(
	\begin{array}{cr}
		1&\\
		&;\frac{-S_0 xt}{1-S_0 t}\\
		\alpha+1&
	\end{array}\right)\nonumber \\
	&=&\sum_{n\geq 0}\frac{L_{n}^{(\alpha)}(x)}{(\alpha+1)_{n}}\left(-S_0 t \right)^n.
\end{eqnarray*}
Finally, the GFs in Theorem \ref{thLaguerre} are obtained after the transformation $t\to t/S_0$.
\end{proof}
\subsection{Bessel polynomials}
The general Bessel polynomials $y_{n}^{(\alpha,\beta)}(x)$, with $\alpha \notin \mathbb{Z}^{-}_{0}$ and $\beta \neq 0$, are expressed by \cite{srivastava2023introductory}
 \begin{equation*}
 	y_{n}^{(\alpha,\beta)}(x)=\sum_{k=0}^{n}\left(^{n}_{k}\right) \left(n+\alpha-1\right)_{k}\left(\frac{x}{\beta} \right)^{k}=\sum_{k=0}^{n}\frac{\left(-n\right)_{k}}{k!} \left(n+\alpha-1\right)_{k}\left(\frac{-x}{\beta} \right)^{k} 
 \end{equation*}
where $y_n(x)=	y_{n}^{(2,2)}(x)=y_{n}^{(2,1)}(x/2)$ are Bessel polynomials. As $\beta$ is a scaling parameter, we will consider the monic polynomials $$y_{n}^{(\alpha)}(x):=2^{n}y_{n}^{(\alpha,1)}(x/2)/(n+\alpha-1)_{n}.$$ 
\begin{theorem}
	The only generating function of the form \eqref{gf0} of the monic Bessel polynomials $y_{n}^{(\alpha)}(x)$ is
	\begin{eqnarray}
	&&
 \sum_{n\geq 0}\left(\frac{3}{2}\right)_ny_{n}^{(3)}(x)\left( 2t\right) ^n	\cong	{\left(1-t\right)^{-2}}\, {}_2F_0\left(
		\begin{array}{cr}
			\frac{3}{2}\,\,\,\,1&\\
			&;\frac{2 xt}{\left(1- t\right)^{2}}\\
			-&
		\end{array}\right),\label{B1}\\
	&&	\mbox{with }A(t)=(1- t)^{-2}, \,R(t)=0\mbox{  and }F(t)={}_2F_0\left(
	\begin{array}{cr}
		\frac{3}{2}\,\,\,\,1&\\
		&;2t\\
		-&
	\end{array}\right).\nonumber
	\end{eqnarray}
	\end{theorem}
\begin{proof}
	The Jacobi-Szeg\"o parameters of ${y}_{n}^{(\alpha)}(x)$ are given by
\begin{equation*}
	\beta_{n}={\frac {2(2-\alpha)}{ \left( 2\,n+\alpha \right)  \left( 2\,n+\alpha-2
			\right) }}
			\end{equation*}
	 and
\begin{equation*}	 
	   \omega_{n}=-{\frac {4n \left( n+\alpha-2 \right) }{ \left( 2\,n+\alpha-1 \right) 
	 		\left( 2\,n+\alpha-3 \right)  \left( 2\,n+\alpha-2 \right) ^{2}}}
\end{equation*}
with $\beta_0=-2/\alpha$ and $\omega_0=0$.\\
From \eqref{eqi} we have for $n\geq 0$
\begin{equation}\label{eqBes1}
	A^{0}_{n+1}=-S_0\frac{\alpha_{n}}{\alpha_{n+1}}+R_0=\frac{1}{n+1}\sum_{k=0}^{n}\beta_k=\frac{-2}{2n+\alpha}.
\end{equation}
Similarly we can show that $A_{n}^{k}=0$ for $n\geq k+1\geq 3$ or equivalently
\begin{eqnarray}\label{eqBes2}
&&-{\frac {R_{{0}}S_{{k}}}{S_{{0}}}}+R_{{k}}-{\frac {S_{{k-2}}S_{{0
		}} \left( R_{{0}} \left(\alpha-2 \right) +2 \right)  \left( R_{{0}}
		\left(\alpha-2 \right) -2 \right) R_{{0}}}{2 \left( R_{{0}}-2
		\right)  \left( R_{{0}}+2 \right)  \left( 2\,R_{{0}}n-2\,R_{{0}}k+R_{
			{0}}\alpha+2 \right) }}\\
		&&+{\frac { \left(\alpha-2 \right) S_{{k-1}}S_{
			{0}}-2\,S_{{k}}}{ \left( 2\,n-2\,k-2+\alpha \right) S_{{0}}}}+{
	\frac { \left(\alpha-2\right)  \left(\left(\alpha-2\right) S_{
			{k-2}}S_{{0}}-2\,S_{{k-1}} \right) }{2(2\,n-2\,k+\alpha)}}\nonumber\\
		&&-{\frac {S_{{k-2}}S_{{0}} \left( \alpha-1 \right)  \left( \alpha-3 \right) }{
		2\left( 2\,n-2\,k-1+\alpha \right)  \left( R_{{0}}+2 \right) }}+{
	\frac {S_{{k-2}}S_{{0}} \left( \alpha-1 \right)  \left( \alpha-3
		\right)}{ 2\left( 2\,n-2\,k+1+\alpha \right)  \left( R_{{0}}-2
		\right) }}=0.\nonumber
\end{eqnarray}
Remark that from \eqref{eqBes1}, $S_{0}\neq 0$, and \eqref{eqBes2} gives us $\alpha=3$ or $\alpha=1$, $R_0=0$, $R_k=0$ and $S_{k-1}=(\alpha-2)S_0S_{k-2}/2$ for $k\geq 2$. Equation \eqref{eqii} gives $R_1=(3-\alpha)S_0$ and \eqref{eqiii} $\alpha=3$. Equation \eqref{eqiv} is verified, so $S_{k}=S_0 \left( S_{0}/2\right)^k$ and we get
$A(t)=(1-S_0 t/2)^{-2}$, $R(t)=0$ and from \eqref{eqBes1} we have
\begin{equation*}
\alpha_n=S_0^n\left(\frac{3}{2}\right)_n, \;n\geq 0
\end{equation*}
which leads to 
\begin{eqnarray*}
	F(t)={}_2F_0\left(
	\begin{array}{cr}
		\frac{3}{2}\,\,\,\,1&\\
		&;S_0 t\\
		-&
	\end{array}\right),
\end{eqnarray*}
and then
\begin{eqnarray*}
A(t)F(xtA(t)-R(t))&=&{\left(1-\frac{S_0 t}{2}\right)^{-2}}\, {}_2F_0\left(
	\begin{array}{cr}
		\frac{3}{2}\,\,\,\,1&\\
		&;\frac{S_0 xt}{\left(1- \frac{S_0t}{2}\right)^{2}}\\
		-&
	\end{array}\right)\nonumber\\
	&\cong& \sum_{n\geq 0}\left(\frac{3}{2}\right)_ny_{n}^{(3)}(x)\left( S_0t\right)^n .
\end{eqnarray*}
The GF for this case with $S_0=2$ can be found in \cite[eq (6) page 294 with $c=2$]{rainville1971special}.

\end{proof}

\section{ Obtaining GFs by derivation}\label{sec4}
We can apply the derivative of order $m$, $\frac{d^{m}}{dx^{m}}$, to the GFs of previous sections to get other GFs of Rainville type \eqref{gf1}. Effectively, by applying $\frac{d^{m}}{dx^{m}}$ to \eqref{gf0} we get
\begin{equation}\label{GFderive}
	A^{m+1}(t)F^{(m)}(xtA(t)-R(t))=\sum_{n\geq 0}\alpha_{n+m}\frac{d^{m}}{dx^{m}}P_{n+m}(x)t^n
\end{equation}
where $B(t)=A^{m+1}(t)$ and  $F^{(m)}(z):=\frac{d^{m}}{dz^{m}}F(z)$.\\
The $m$ derivatives of the classical monic orthogonal polynomials are given by
\begin{eqnarray*}
&&	\frac{d^{m}}{dx^{m}}P_{n}^{\left( \alpha,\beta\right) }(x)=\frac{n!}{\left( n-m\right)!}P_{n-m}^{\left( \alpha+m,\beta+m\right) }(x)\\
&&	\frac{d^{m}}{dx^{m}}C_{n}^{\lambda}(x)=\frac{n!}{\left( n-m\right)!}C_{n-m}^{\lambda+m}(x)\\
&&	\frac{d^{m}}{dx^{m}}U_{n}(x)=\frac{d^{m}}{dx^{m}}C_{n}^{1}(x)=\frac{n!}{\left( n-m\right)!}C_{n-m}^{1+m}(x)\\
&&	\frac{d^{m}}{dx^{m}}H_{n}(x)=\frac{n!}{\left( n-m\right)!}H_{n-m}(x)\\
&&	\frac{d^{m}}{dx^{m}}L_{n}^{\alpha}(x)=\frac{n!}{\left( n-m\right)!}L_{n-m}^{\alpha+m}(x)\\
&&	\frac{d^{m}}{dx^{m}}y_{n}^{\alpha}(x)=\frac{n!}{\left( n-m\right)!}y_{n-m}^{\alpha+2m}(x).
\end{eqnarray*}
The interesting cases concerns the GFs \eqref{U1} for the Chebyshev polynomials of the second kind and \eqref{H2} for Hermite polynomials. The corresponding results are given in the two theorems below.

\begin{theorem}\label{thUG} The following GF holds true
	\begin{eqnarray}\label{U22}
		&&\sum_{n=0}^{+\infty}\frac{(n+m)!2^{n+m}}{n!m!}[U_{n+m}\left(\rho\right)+(\tfrac{\alpha_1}{2}-\rho)U_{n+m-1}\left(\rho\right)] C^{m+1}_n(x)t^n\nonumber\\
		&&\qquad=\frac{c_1}{\left( t_1-xt-\left( \frac{t_1}{4}-\frac{\rho}{2}\right) t^2\right)^{m+1} }+\frac{c_2}{\left(t_2-xt-\left( \frac{t_2}{4}-\frac{\rho}{2}\right) t^2\right)^{m+1}}.\;\;\;\;\;\;\;\;
	\end{eqnarray}
	where 
	\begin{equation}\label{t1t2}
		t_1=\rho-\sqrt{\rho^2-1},\; t_2=\rho+\sqrt{\rho^2-1},\; 
		\end{equation}
		and 
		\begin{equation}\label{c1c2}
		 c_1=\frac{1-\left(\alpha_1 -2\rho \right)t_1}{t_2-t_1},\; c_2=\frac{1-\left(\alpha_1 -2\rho \right)t_2}{t_1-t_2}.
	\end{equation}
\end{theorem}
\begin{proof} For the GF \eqref{U1}  we have 
	\begin{equation*}
		F ( t ) ={\frac { 1+\left(\alpha_1 -2\rho \right) t
			}{1-2\rho t+{t}^{2}}},\;\; A(t)=\frac{4}{4-t^2} \text{ and } R(t)= \frac{2\rho t^2}{4-t^2}.
	\end{equation*}
	To calculate the $m$ derivative of $F(t)$ we write $F(t)$ in its partial fraction decomposition
	\begin{equation*}
		F (t) ={\frac { 1+\left(\alpha_1 -2\rho \right) t
			}{1-2\rho t+{t}^{2}}}=\frac{c_1}{t_1-t}+\frac{c_2}{t_2-t}
	\end{equation*}
	where $t_1$, $t_2$, $c_1$ and $c_2$ are given in \eqref{t1t2} and \eqref{c1c2}.  \\	
	Then
	\begin{equation*}
		F^{(m)} ( t ) =\frac{c_1 m!}{\left( t_1-t\right)^{m+1} }+\frac{c_2 m!}{\left( t_2-t\right)^{m+1}}
	\end{equation*}
	and \eqref{GFderive} becomes
	\begin{eqnarray*}
		&&\sum_{n=0}^{+\infty}\frac{(n+m)!}{n!}\alpha_{n+m} C^{m+1}_n(x)t^n\\
		&&\qquad\qquad\qquad= m! A(t)^{m+1}\left( \tfrac{c_1}{\left(t_1-xtA(t)+R(t)\right)^{m+1} }+\tfrac{c_2}{\left( t_2-xtA(t)+R(t)\right)^{m+1}}\right)\;\;\;\;\;
	\end{eqnarray*}
	leading to the GF \eqref{U22}.
\end{proof}
\begin{theorem} The following GFs hold true:
	\begin{eqnarray}
		&&* \sum_{n=0}^{+\infty}\frac{2^n}{n!}H_{n+m}(\rho) H_n(x)t^n=\left( 1-t^2\right)^{-\frac{m+1}{2}} e^{{\rho}^{2}-\frac{(\rho-xt)^2}{1-t^2}} H_{m}\left( \frac{\rho-xt}{\sqrt{1-t^2}}\right) \;\;\;\;\;\;\label{H30}\\
		&&* \sum_{n=0}^{+\infty}\frac{2^n}{n!}H_{n+m}(\rho) t^n=H_{m}\left(\rho-t\right)  e^{2\rho t-t^2} \;\;\;\;\;\;\label{H31}\\
		&&* \sum_{n=0}^{+\infty}\frac{2^n}{n!}H_{n+m-1}(\rho,1) H_n(x)t^n=\left( 1-t^2\right)^{-\frac{m+1}{2}}\left[H_{m-1}\left(\frac{\rho-xt}{\sqrt{1-t^2}},1\right)\right.\nonumber\\
		&&\left.\qquad\qquad+\sqrt{\pi} \left( {\rm erfi} \left(\rho\right) -{\rm erfi} \left(\frac{\rho-xt}{\sqrt{1-t^2}}\right) \right) e^{-\frac{(\rho-xt)^2}{1-t^2}} H_{m}\left( \frac{\rho-xt}{\sqrt{1-t^2}}\right) \right]\;\;\;\;\;\;\label{H32}\\
		&&* \sum_{n=0}^{+\infty}\frac{2^n}{n!}H_{n+m-1}(\rho,1) t^n=H_{m-1}\left(\rho-t,1\right)\nonumber\\
		&&\qquad\qquad\qquad\qquad\qquad+\sqrt{\pi}\left( {\rm erfi} \left(\rho\right) -{\rm erfi} \left(\rho-t\right) \right) e^{-(\rho-t)^2} H_{m}\left( \rho-t\right). \;\;\;\;\;\;\;\;\label{H33}
	\end{eqnarray}
\end{theorem}
\begin{proof} We apply $\frac{d^{m}}{dx^{m}}$ to \eqref{H2} to get
	\begin{eqnarray}\label{H3}
		&&\sum_{n=0}^{+\infty}\frac{(n+m)!}{n!}\alpha_{n+m} H_n(x)t^n\nonumber\\
		&&\qquad\qquad\qquad=\left( 1-t^2\right)^{-\frac{m+1}{2}} F^{(m)}\left(\frac{xt}{\sqrt{1-t^2}}-\rho\left( \frac{1}{\sqrt{1-t^2}}-1\right)\right) \;\;\;\;\;\;
	\end{eqnarray}
	where
	\begin{eqnarray*}\label{H}
		F(t)&=&\left(1+ \tfrac{\alpha_1-2\rho}{2}\sqrt{\pi}e^{-{{\rho}^{2}}}   \left[{\rm erfi} \left(\rho\right) -{\rm erfi} \left(\rho-t\right) \right]\right)e^{t(2\rho-t)}\\\;\;\;
		&&=\left(\rho_2- \rho_1{\rm erfi} \left(\rho-t\right)\right)e^{-(\rho-t)^2}\;\;\;
	\end{eqnarray*}
	with $\rho_1=\tfrac{\alpha_1-2\rho}{2}\sqrt{\pi}$ and $\rho_2=e^{{{\rho}^{2}}}+ \rho_1 {\rm erfi} \left(\rho\right)$.
	
	-- If $\alpha_1=2\rho$, i.e. $\rho_1=0$, then $F(t)=e^{{{\rho}^{2}}}e^{-(\rho-t)^2}$ and 
	\begin{equation*}
		F^{(m)}(t)=e^{{{\rho}^{2}}}\frac{d^{m}}{dt^{m}}\left[  e^{-(\rho-t)^2}\right]=e^{{{\rho}^{2}}}(-1)^{m}\frac{d^{m}}{dz^{m}}\left[  e^{-z^2}\right]_{z=\rho-t}
	\end{equation*}
	where clearly for $t=0$ the GF \eqref{H3} gives us
	\begin{equation*}
		m!\alpha_{m}=2^{m}H_{m}(\rho)=(-1)^{m}e^{{{\rho}^{2}}}\frac{d^{m}}{d\rho^{m}} e^{-\rho^2}.
	\end{equation*}
	So,  $F^{(m)}(t)=2^{m}H_{m}(\rho-t)F(t)$ leading to the GF \eqref{H30}. To obtain the known GF \eqref{H31} \cite[p. 197]{rainville1971special}, we make the change $x\to x/t$ and then the limit $t\to 0$ in the GF \eqref{H30}.
	
	-- If $\alpha_1\neq 2\rho$, we use the identity
	\begin{equation*}
		\frac{d^{m}}{dz^{m}}{\rm erfi} \left(z\right)=\frac{2}{\sqrt{\pi}}\frac{d^{m-1}}{dz^{m-1}}e^{z^2}=\frac{2^{m}(-i)^{m-1}}{\sqrt{\pi}}H_{m-1}\left(iz\right)e^{z^2},\;m\geq 1
	\end{equation*}
	and Leibniz rule to get
	\begin{eqnarray*}\label{H4}
		F^{(m)}(t)&=&2^{m}F(t)	H_{m}\left(\rho-t\right)+\sum_{k=1}^{m}\tfrac{m!}{k!(m-k)!}\left(- \rho_1\tfrac{d^{k}}{dt^{k}}{\rm erfi} \left(\rho-t\right) \right)\tfrac{d^{m-k}}{dt^{m-k}} e^{-(\rho-t)^2}\;\;\;\;\;\;\;\;\;\;\;\;\\
		&=&2^{m}F(t)H_{m}\left(\rho-t\right)+\frac{2^{m}\rho_1}{\sqrt{\pi}}\sum_{k=1}^{m}\tfrac{i^{k-1}m!}{k!(m-k)!}H_{k-1}\left(i\left( \rho-t\right)\right)H_{m-k}\left(\rho-t\right)\;\;\;\;\;\;\;\;\;\\
		&=&2^{m}F(t)H_{m}\left(\rho-t\right)+2^{m-1}\left( \alpha_{{1}}-2\rho\right)H_{m-1}(\rho-t,1)
	\end{eqnarray*}
	where we have used the following expression, see \cite[formula 4.15]{AskeyWimp1984}, 
	\begin{eqnarray}
		H_{m-1}(\rho,1)=\sum_{k=1}^{m}\tfrac{i^{k-1}m!}{k!(m-k)!}H_{k-1}\left(i\rho\right)H_{m-k}\left(\rho\right).
	\end{eqnarray}
	In this case we obtain the GF \eqref{H32} where for $m=0$ we retrieve the GF \eqref{H44} and by the change $x\to x/t$ and the limit $t\to 0$ we get, after replacing $x$ by $t$, the GF \eqref{H33} for the associated Hermite polynomials which is similar to \eqref{H31}.
\end{proof}
\begin{remark}
	Note that the GF \eqref{H30} is a special case of certain extensions of Mehler's formula \cite{Carlitz1970Mehler,tapia2024generalization}. Also, the GF \eqref{U2} as a sub-case of \eqref{U22} can be found in \cite{CARLITZL1972} and \cite[Eq. 2.10]{szablowski2022multivariate}. Whereas, the GFs \eqref{H44}, \eqref{H32}, \eqref{H33} and \eqref{U22} seems to be new in our knowledge.
\end{remark}
Now for the other cases we can easily check that the GFs \eqref{H1}, \eqref{L2}, \eqref{J1}, \eqref{Ultra1} and \eqref{Ultra2} remain the same and that for the GFs \eqref{T1} and \eqref{U3} we obtain \eqref{Ultra1}. Concerning the GFs \eqref{J2}, \eqref{J3}, \eqref{L1} and \eqref{B1}; the application of the $m$ derivative produces known GFs, but for sake of completeness we list them below.

$\bullet$ To the GF \eqref{J2}:
\begin{eqnarray*}\label{J22}
	&& \left( 1+\frac{t}{2}\right)^{-2-2m} {}_2F_1\left(
	\begin{array}{cr}
		\frac{3}{2}+m\,\,\,\,1+m&\\
		&;\frac{(x+1)t}{\left( 1+\frac{t}{2}\right)^{2}}\\
		2+m&
	\end{array}\right)\nonumber\\
	&&\qquad\qquad\qquad=\sum_{n\geq0}\frac{\left(\frac{3}{2}+m\right)_{n}\left(1+m\right)_{n}}{(2+m)_{n}n!}\tilde{P}_{n}^{(m,m+1)}(x)t^{n}.\;\;\;\;\;\;\;
\end{eqnarray*}

$\bullet$ To the GF \eqref{J3}:
\begin{eqnarray*}\label{J33}
	&&{\frac {1}{ \left( 1-t\right) ^{2+2m}} }{}_2F_1\left(
	\begin{array}{cr}
		\frac{3}{2}+m\,\,\,\,1+m&\\
		&;{\frac {2t
				\left( x-R_{{0}} \right) }{ R_{{0}}\left(1-t\right) ^{2}}}
		\\
		\frac{3}{2}+\frac{2\alpha-1}{2R_0}+m&
	\end{array}\right)\\
	&&\qquad\qquad\qquad=\sum_{n\geq 0}\tfrac{\left(\frac{3}{2}+m\right)_n \left( 1+m\right)_n }{\left(\frac{3}{2}+\frac{2\alpha-1}{2R_0}+m\right)_{n}n! }\tilde{P}_{n}^{(\alpha+m,1-\alpha+m)}(x)\left(\tfrac{2t}{R_0}\right)^n\;\;\;
\end{eqnarray*}
where by the substitution $\alpha+m\to \alpha$, $1-\alpha+m \to \beta$, $x\to x+1$ and $R_0=\rho=1$ we recover the GF obtained in \cite[Eq 2.14]{ismail1974obtaining},
\begin{eqnarray}\label{J333}
	&&{\frac {1}{ \left( 1-t\right) ^{\alpha+\beta+1}} }{}_2F_1\left(
	\begin{array}{cr}
		\frac{\alpha+\beta+2}{2}\,\,\,\,\frac{\alpha+\beta+1}{2}&\\
		&;{\frac {2t
				x }{\left(1-t\right) ^{2}}}
		\\
		\alpha+1&
	\end{array}\right)\nonumber\\
	&&\qquad\qquad\qquad=\sum_{n\geq 0}\tfrac{\left(\frac{\alpha+\beta+2}{2}\right)_n \left( \frac{\alpha+\beta+1}{2}\right)_n }{2^{-n}\left(\alpha+1\right)_{n}n! }\tilde{P}_{n}^{(\alpha,\beta)}(x+1) t^n.\;\;\;
\end{eqnarray}
We note that \eqref{J333} is not a GF of Boas-Buck type as stated in \cite[Eq 2.14]{ismail1974obtaining}, but of Rainville type. Also by choosing $\alpha=0$, $R_0=-1$ and $t\to -t/2$ we get \eqref{J22}.

$\bullet$ To the GF \eqref{L1}:
\begin{eqnarray*}\label{L11}
	\left( 1- t\right) ^{-1-m}{}_1F_1\left(
	\begin{array}{cr}
		1+m&\\
		&;\frac{-xt}{ 1-t}\\
		\alpha+m+1&
	\end{array}\right)
	&=&\sum_{n\geq 0}\frac{\left( 1+m\right)_{n} (-1)^n}{(\alpha+m+1)_{n}n!}L_{n}^{(\alpha+m)}(x)t^n.\;\;\;\;\;\;\;\;\;\;\;\;
\end{eqnarray*}
where by the substitution $\alpha+m\to \alpha$ and $1+m \to c$ we recognize equation 2.18 in \cite{ismail1974obtaining}.

$\bullet$ To the GF \eqref{B1}:
\begin{eqnarray*}\label{B11}
	&&{\left(1-t\right)^{-2-2m}}\, {}_2F_0\left(
	\begin{array}{cr}
		\frac{3}{2}+m\,\,\,\,1+m&\\
		&;\frac{2 xt}{\left(1- t\right)^{2}}\\
		-&
	\end{array}\right)
	\\
	&&\qquad\qquad\qquad\cong \sum_{n\geq 0}\frac{\left(\frac{3}{2}+m\right)_n\left(1+m\right)_n }{n!}y_{n}^{(3+2m)}(x)\left( 2t\right) ^n\;\;\;\;\;\;\;
\end{eqnarray*}
which can be compared to equation 2.20 in \cite{ismail1974obtaining} by the change $(2+2m) \to c$.

\bibliographystyle{spmpsci}
\bibliography{mesk_zahaf}

\end{document}